\documentclass[12pt]{elsarticle}
\usepackage{amsfonts,enumerate,amsmath,amssymb,mathrsfs,amsbsy}
\usepackage[colorlinks=true]{hyperref}

\usepackage{caption}
\biboptions{numbers,sort&compress}

\usepackage{xcolor}

\def\qed{\strut\hfill $\Box$}
\newtheorem{thm}{Theorem}[section]

\newtheorem{lem}[thm]{Lemma}

\newtheorem{rem}[thm]{Remark}
\newtheorem{defn}[thm]{Definition}

\newcommand{\thmref}[1]{Theorem~{\rm \ref{#1}}}
\newcommand{\lemref}[1]{Lemma~{\rm \ref{#1}}}

\def\para#1{\vskip .4\baselineskip\noindent{\bf #1}}
\numberwithin{equation}{section}
\allowdisplaybreaks \allowdisplaybreaks[4]
\begin{document}
\begin{frontmatter}	
	\title{Strong averaging principles for a class of non-autonomous slow-fast systems of SPDEs with polynomial growth}
	
	\author[mymainaddress]{Ruifang Wang}
	\ead{wrfjy@yahoo.com}
	 
	\author[mymainaddress,myfivearyaddress]{Yong Xu\corref{mycorrespondingauthor}}
	\cortext[mycorrespondingauthor]{Corresponding author}
	\ead{hsux3@nwpu.edu.cn}
	
	
	
	
	\address[mymainaddress]{Department of Applied Mathematics, Northwestern Polytechnical University, Xi'an, 710072, China}
	\address[myfivearyaddress]{MIIT Key Laboratory of Dynamics and Control of Complex Systems, Northwestern Polytechnical University, Xi'an, 710072, China	}
	

	%

	\begin{abstract}
	In this work, we study a class of non-autonomous two-time-scale stochastic reaction-diffusion equations driven by Poisson random measures, in which the coefficients satisfy the polynomial growth condition and local Lipschitz condition. First, the existence and uniqueness of the mild solution are proved by constructing auxiliary equations and using the technique of stopping time. Then, consider that the time dependent of the coefficients, the averaged equation is redefined by studying the existence of time-dependent evolution family of measures associated with the frozen fast equation. Further, the slow component strongly converges to the solution of the corresponding averaged equation is verified by using the classical Khasminskii method.
	\end{abstract}		
\end{frontmatter}

\section{Introduction}\label{sec-1}
In this paper, we are concerned with the following  non-autonomous slow-fast systems
of stochastic partial differential equations (SPDEs) on a bounded domain $ \mathcal{O} $  of  $ \mathbb{R}^d\left( d\ge 1 \right)  $: 
\small\begin{eqnarray}\label{orginal1}
\begin{aligned}
\begin{cases}
\frac{\partial X^{\epsilon}}{\partial t}\left( t,\xi \right) &=\mathcal{A}_1  X^{\epsilon}\left( t,\xi \right) +b_1\left(  \xi ,X^{\epsilon}\left( t,\xi \right),Y^{\epsilon}\left( t,\xi \right) \right) +f_1\left(  \xi ,X^{\epsilon}\left( t,\xi \right) \right) \frac{\partial W ^{Q_1}}{\partial t}\left( t,\xi \right) \\
&\quad+\int_{\mathbb{Z}}{g_1\left( \xi ,X^{\epsilon}\left( t,\xi \right) ,z \right) \frac{\partial \tilde{N}_1}{\partial t}\left( t,\xi ,dz \right)},\\
\frac{\partial Y^{\epsilon}}{\partial t}\left( t,\xi \right) 
&= \frac{1}{\epsilon}\left[ \left( \mathcal{A}_2\left( t \right) -\alpha \right) Y^{\epsilon}\left( t,\xi \right) +b_2\left( t,\xi ,X^{\epsilon}\left( t,\xi \right) ,Y^{\epsilon}\left( t,\xi \right) \right) \right] \\
&\quad+\frac{1}{\sqrt{\epsilon}}f_2\left( t, \xi, Y^{\epsilon}\left( t,\xi \right) \right) \frac{\partial W ^{Q_2}}{\partial t}\left( t,\xi \right) +\int_{\mathbb{Z}}{g_2}\left( t, \xi, Y^{\epsilon}\left( t,\xi \right) ,z \right) \frac{\partial \tilde{N}_{2}^{\epsilon}}{\partial t}\left( t,\xi ,dz \right),\\
X^{\epsilon}\left( 0,\xi \right) &=x\left( \xi \right), \quad 
Y^{\epsilon}\left( 0,\xi \right) =y\left( \xi \right), \quad
\xi \in \mathcal{O},\cr
\mathcal{N}_1X^{\epsilon}\left( t,\xi \right) &= \mathcal{N}_2Y^{\epsilon}\left( t,\xi \right) =0, \quad t\ge 0, \quad \xi \in \partial \mathcal{O},
\end{cases}
\end{aligned}
\end{eqnarray}
where $ W ^{Q_1}, W ^{Q_2}$ and $\tilde{N}_1, \tilde{N}_{2}^{\epsilon} $ are mutually independent Wiener processes and Poisson random measures,  $ 0<\epsilon \ll 1 $ is a  small parameter and $ \alpha $  is a sufficiently large fixed constant. In addition, $\mathcal{N}_i (i=1,2)$ are the boundary operators, which can be either the identity operator (Dirichlet boundary condition) or the first order operator (coefficients satisfying a uniform nontangentiality condition). The stochastic perturbations of the equations define on the same complete stochastic basis $\big( \varOmega ,\mathcal{F},\left\lbrace \mathcal{F}_t\right\rbrace _{t\geq0},\mathbb{P} \big) $  and the specific introduction will be given in Section \ref{sec-2}.

The averaging principle is an effective method to analysis the slow-fast systems, which can simplify the system by constructing the averaged equation. In 1961, Bogolyubov and  Mitropolskii \cite{ Bogolyubov1961asymptotic}  studied the averaging principle, giving the first rigorous result for the deterministic case.  
Then, Khasminskii \cite{khas1968on} established the averaging principle for stochastic differential equations (SDEs) in 1968. Since then, the averaging principle became an active area of research. Givon \cite{givon2007strong}, Freidlin and Wentzell \cite{freidlin2012random},  Duan \cite{duan2014effective}, Xu and his co-workers \cite {xu2011averaging, xu2015approximation, xu2017stochastic} studied the averaging principle of SDEs. In recent years, many scholars also investigated the averaging principle of SPDEs. Such as: Stochastic heat equation has been researched by Cerrai and Freidlin \cite{Cerrai2009Averaging}, Cerrai \cite{cerrai2009khasminskii,cerrai2011averaging}, Wang and Roberts \cite{wang2012average}, Fu and co-workers \cite{Fu2011AN}, Bao and co-workers\cite{Bao2017Two}, et al. Moreover, Dong and co-workers \cite{dong2018averaging} studied the Stochastic Burgers equations. Fu and co-workers \cite{Fu2015Strong}, Pei and co-workers \cite{pei2017two} researched the Stochastic wave equation. Gao studied the Stochastic Klein-Gordon equation \cite{gao2019averaging} and Stochastic Schr\"{o}dinger equation \cite{gao2018averaging}.
The above studies about the averaging principle are based on autonomous systems,  but in practical problems, the parameters of the system often depend on time. Non-autonomous system can depict some actual models better, which has made itself attract more and more attention of scholars.
In 2017, Cerrai \cite{cerrai2017averaging} studied the averaging principle for non-autonomous slow-fast systems of stochastic reaction-diffusion equations driven by Gaussian noises. Then, Liu \cite{liu2020averaging} researched the averaging principle for non-autonomous stochastic differential equations driven by Gaussian noises and Xu \cite{Xu2018Averaging} studied the averaging principle for non-autonomous slow-fast systems of stochastic reaction-diffusion equations driven by Gaussian noises and Poisson random measures.


In this work, we also concerned with a non-autonomous slow-fast systems of reaction-diffusion equations
driven by Gaussian noises and Poisson random measures. 
Consider that the Lipschitz condition and linear growth condition are too strict to study the validity of the averaging principle in many other relevant cases.
For example, one of the reaction-diffusion equations is the Fitzhugh-Nagumo or Ginzburg-Landau type, in which the coefficients of these systems satisfy  the polynomial growth condition, has appeared in the fields of biology and physics and attracted considerable attention.
As a consequence, comparing with the work of Xu \cite{Xu2018Averaging}, the coefficients of the original equation (\ref{orginal1}) in this work are extended to satisfy the polynomial growth condition and the local Lipschitz condition. In order to deal with it, subdifferential is considered to obtain some a-priori bounds 
and the stopping technique will be used very frequently. 
Moreover, the existence of time-dependent evolution family of measures associated with the frozen fast equation is also studied to define the averaged equation and the technique of time discretization and truncation is used to obtain the averaging principle. 

The purpose of this work is to prove that the slow component strongly converges to the solution of the corresponding averaged equation under some reasonable assumptions. 
That is, for any  $ T>0 $ and  $ p\geq 1 $, we have 
\begin{align} 
\underset{\epsilon \rightarrow 0}{\lim}\mathbb{E}\Big( \underset{t\in \left[ 0,T  \right]}{{\sup}}\left\|  X^{\epsilon}\left( t \right) -\bar{X}\left( t \right) \right\| ^p\Big) =0,
\end{align}
where  $ \bar{X} $ is the solution of the corresponding averaged equation:
\begin{eqnarray}\label{en15}
\begin{split}
\begin{cases}
\frac{\partial \bar{X}}{\partial t}\left( t,\xi \right) &=\mathcal{A}_1  \bar{X}\left( t,\xi \right) +\bar{B}_1\left( \bar{X}\left( t \right) \right)\left( \xi\right)  +f_1\left(  \xi ,\bar{X}\left( t,\xi \right) \right) \frac{\partial W ^{Q_1}}{\partial t}\left( t,\xi \right) \\
&\quad+\int_{\mathbb{Z}}{g_1\left(  \xi ,\bar{X}\left( t,\xi \right) ,z \right) \frac{\partial \tilde{N}_1}{\partial t}\left( t,\xi ,dz \right)},\\
\bar{X}\left( 0,\xi \right) &=x\left( \xi \right), \quad \xi\in\mathcal{O}, \quad
\mathcal{N}_1\bar{X}\left( t,\xi \right) =0, \quad t\ge 0,\quad \xi \in \partial\mathcal{O},
\end{cases} 
\end{split}
\end{eqnarray} 
with 
\begin{eqnarray}\label{en18}
\bar{B}_1\left( x \right) :=\underset{T\rightarrow \infty}{\lim}\frac{1}{T}\int_0^T{\int_{L^2\left(\mathcal{O}\right)}{B_1\left( t,x,y \right)}\mu _{t}^{x}\left( dy \right) dt, \quad x\in L^2\left(\mathcal{O}\right)},
\end{eqnarray}
and $\left\{ \mu _{t}^x  \right\}_{t \in \mathbb{R} }$ is an evolution family of probability measures for the fast equation with frozen slow component $ x $.

This work consists of the following contents. In Section \ref{sec-2}, we introduce some notations, assumptions and presents the main results of this work. Then, the existence and uniqueness of  solutions is proved and some a-priori estimates for the solution of  (\ref{orginal1}) are given in Section \ref{sec-3}. 
Finally, we define the averaged equation by studying the frozen fast equation and give
the detailed proof of our main
result by using the generalized Khasminskii method  in Section \ref{sec-4}.
In this paper, $ C>0 $ with or without subscripts represents a general constant, the value of which may vary for different cases.

\section{Notations, assumptions and preliminaries}\label{sec-2}

Denote $  E  $ is a Banach space  and endowed with the following sup-norm 
$$
\left\| x \right\|  =\underset{\xi \in \bar{\mathcal{O}}}{\text{sup}}\left| x\left( \xi \right) \right|,
$$
and the duality $ \left< \cdot ,\cdot \right> $. The norm of the product space $   E \times E $ denote as   
$$
\left\| x \right\| _{  E \times   E }=\big( \left\| x_1 \right\|  ^{2}+\left\| x_2 \right\|  ^{2} \big) ^{\frac{1}{2}},
$$
and the corresponding duality of the product space $   E \times E $ is $ \left< \cdot ,\cdot \right> _{  E \times  E } $. 


Let $  X  $ be any 	Banach space, denote $\mathcal{L}\left(  X  \right)$ is the space of the bounded linear operators in $ X, $  and denote $\mathcal{L}_2\left( X \right)$ the subspace of Hilbert-Schmidt operators, endowed with the norm $$
\left\| Q \right\| _2=\sqrt{Tr\left[ Q^{\ast}Q \right]}.
$$ 
Moreover, $  \mathbb{D}\left( \left[ s,T \right] ;  X  \right) $ denotes the space of all c\`adl\`ag path from $ \left[ s,T \right] $ into $  X  $.


Next, we introduce some notations about subdifferential, which will be used in what follows (reader can see 
\cite[Appendix D]{prato2014stochastic} for all definitions and details).  
The subdifferential $ \partial \left\| x \right\| $ of $ \left\| x \right\|   $ is defined as 
$$
\partial \left\| x \right\|  :=\left\{ h\in   E ^*;\left\| h \right\| _{  E ^*}=1,\left< x,h \right>  =\left\| x \right\|   \right\},
$$
where $  E ^* $ is the dual space of $ E $. 

For any $ x \in E, $ we denote
\begin{align}
M_{x}=\left\lbrace  \xi \in \bar{\mathcal{O}}:|x(\xi)|=\left\| x\right\| \right\rbrace .
\nonumber
\end{align}
Then, for any $ x\in E \backslash \left\lbrace 0\right\rbrace ,$ we denote
\begin{align}
\mathcal{M}_{x}=\left\{\delta_{x, \xi} \in E^{*} ; \xi \in M_{x}\right\},\nonumber
\end{align}
where
\begin{align}
\left\langle\delta_{x, \xi}, y\right\rangle =\frac{x(\xi) y(\xi)}{\left\| x\right\| }, \quad y \in E,\nonumber
\end{align}
and for $ x = 0, $ we set
\begin{align}
\mathcal{M}_{0}=\left\{h \in E^{*}:\left\| h\right\| _{E^{*}}=1\right\}.\nonumber
\end{align}
It is easy to know that $ \mathcal{M}_{x} \subseteq \partial \left\| x \right\|  $ for every $ x\in E. $ Moreover,
due to the characterization of the subdifferential \cite[Appendix D]{prato2014stochastic}, if $ X:\left[ 0,T\right] \rightarrow E  $  is any differentiable mapping, then
\begin{eqnarray}\label{en01}
\frac{d}{dt}^{-}\left\| X\left( t \right) \right\|  \le \left< X{'}\left( t \right) ,\delta \right>  ,
\end{eqnarray}
for any $ t\in \left[ 0,T\right]  $ and $ \delta \in \mathcal{M}_{ X\left( t \right)} $. 

Analogously, if $ x \in E \times E, $ there also have corresponding notations about subdifferential, and this paper will not give more details.

Now, we assume the space dimension $d>1$, the processes ${\partial W ^{Q_1}}/{\partial t}\left( t,\xi \right)$ and ${\partial W ^{Q_2}}/{\partial t}\left( t,\xi \right)$ in the slow-fast system are the Gaussian noises, and assumed it is white in time and colored in space. Here, $W ^{Q_i}\left( t,\xi \right) \left( i=1,2 \right)$ is the cylindrical Wiener processes, defined as $$
W ^{Q_i}\left( t,\xi \right) =\sum_{k=1}^{\infty}{Q_ie_k\left( \xi \right) \beta _k\left( t \right)}, \quad i=1,2,
$$
where $\left\{ e_k \right\} _{k\in \mathbb{N}}$ is a complete orthonormal basis, $\left\{ \beta _k\left( t \right) \right\} _{k\in \mathbb{N}}$ is a sequence of mutually independent standard Brownian motion defined on the same complete stochastic basis $\big( \varOmega ,\mathcal{F},\left\lbrace \mathcal{F}_t\right\rbrace _{t\geq0},\mathbb{P} \big) $, and $Q_i$ is a bounded linear operator.

Next, we give the definitions of Poisson random measures $\tilde{N}_1\left( dt,dz \right)$ and $\tilde{N}_{2}^{\epsilon}\left( dt,dz \right)$. Let $ \left( \mathbb{Z},\mathcal{B}\left( \mathbb{Z} \right) \right) $  be a given measurable space and  $ Y\left( dz \right)  $ be a  $ \sigma  $-finite measure on it. $ D_{p_{t}^{i}},i=1,2 $ are two countable subsets of $ \mathbb{R}_+ $. Moreover, let
$ p_{t}^{1},t \in  D_{p_{t}^{1}}  $ be a stationary $ \mathcal{F}_t $-adapted Poisson point process on $ \mathbb{Z} $ with characteristic $ v $,
and $ p_{t}^{1},t \in  D_{p_{t}^{2}} $ be the other stationary $ \mathcal{F}_t $-adapted Poisson point process on $ \mathbb{Z} $ with characteristic $ {v}/{\epsilon} $. Denote by $ N_i\left( dt,dz \right),i=1,2 $ the Poisson counting measure associated with $ p_{t}^{i} $, i.e.,
$$
N_i\left( t,\varLambda \right) :=\sum_{s\in D_{p_{t}^{i}},s\le t}{I_{\varLambda}\left( p_{t}^{i} \right)},\quad i=1,2.
$$
Let us denote the two independent compensated Poisson measures
$$\tilde{N}_1\left( dt,dz \right) :=N_1\left( dt,dz \right) - v_1\left( dz \right) dt$$ 
and 
$$\tilde{N}_{2}^{\epsilon}\left( dt,dz \right) :=N_2\left( dt,dz \right) -\frac{1}{\epsilon} v_2\left( dz \right) dt,$$ 
where  $ v_1\left( dz \right) dt$ and $\frac{1}{\epsilon} v_2\left( dz \right) dt$ are the compensators.

In this paper, for any $ t\in \mathbb{R} $, the operators $\mathcal{A}_1  $ and $\mathcal{A}_2  $ independent of $t$ are the second order uniformly elliptic operators  with continuous coefficients on $\bar{\mathcal{O}}$. Moreover, we assume that the operator $ \mathcal{A}_2(t) $ has the following form
\begin{eqnarray}\label{en21}
\begin{split}
\mathcal{A}_2\left( t \right) =\gamma  \left( t \right) \mathcal{A}_2+\mathcal{L} \left( t \right), \quad t\in \mathbb{R}, 
\end{split}
\end{eqnarray}
where the operator $\mathcal{L} \left( t \right) $ is a first order differential operator has the following form
\begin{eqnarray}\label{en22}
\begin{split}
\mathcal{L} \left( t,\xi \right) X\left( \xi \right) =\left< l \left( t,\xi \right) ,\nabla X\left( \xi \right) \right> _{\mathbb{R}^d}, \quad t\in \mathbb{R}, \ \xi \in \bar{\mathcal{O}}.
\end{split}
\end{eqnarray}

Finally, for $ i=1,2 $, denote the realization of the operators $ \mathcal{A}_i $ and $ \mathcal{L}  $  in $ E $  are $ A_i $ and $ L  $, and the operator $ \mathcal{A}_i $ generates an analytic semigroup $ e^{tA_i} $.

Now, we give the following assumptions about the operators $ A_i $ and $ Q_i $ as in \cite{cerrai2017averaging} and \cite{Xu2018Averaging}. 
\begin{enumerate}[({A}1)]
	\item
	\begin{enumerate}
		\item 	The function $\gamma  :\mathbb{R}\rightarrow \mathbb{R}$ is continuous and there exist $\gamma _0, \gamma_1 >0$ such that 
		\begin{eqnarray}\label{en23}
		\begin{split}
		\gamma _0\le \gamma  \left( t \right) \le \gamma_1, \quad t\in \mathbb{R}.
		\end{split}
		\end{eqnarray}
		\item The function $l :\mathbb{R}\times \bar{\mathcal{O}}\rightarrow \mathbb{R}^d$ is continuous and bounded.
	\end{enumerate}	
	\item For $i=1,2$, there exist a complete orthonormal system $\left\{ e_{i,k} \right\} _{k\in \mathbb{N}}$ and two sequences of nonnegative real numbers $\left\{ \alpha _{i,k} \right\} _{k\in \mathbb{N}}$ and $\left\{ \lambda _{i,k} \right\} _{k\in \mathbb{N}}$ such that
	\begin{eqnarray}\label{en24}
	\begin{split}
	A_ie_{i,k}=-\alpha _{i,k}e_{i,k}, \quad Q_ie_{i,k}=\lambda _{i,k}e_{i,k}, \quad k\ge 1,
	\end{split}
	\end{eqnarray}
	and
	\begin{eqnarray}\label{en25}
	\begin{split}
	\kappa _i:=\sum_{k=1}^{\infty}{\lambda _{i,k}^{\rho _i}\left\| e_{i,k} \right\| _{\infty}^{2}}<\infty, \quad \zeta _i:=\sum_{k=1}^{\infty}{\alpha _{i,k}^{-\beta _i}\left\| e_{i,k} \right\| _{\infty}^{2}}<\infty,
	\end{split}
	\end{eqnarray}
	for some constants $\rho _i\in \left( 2,+\infty \right]$ and $\beta _i\in \left( 0,+\infty \right)$ such that
	\begin{eqnarray}\label{en26}
	[{\beta _i\left( \rho _i-2 \right)}]/{\rho _i}<1.
	\end{eqnarray}
	Moreover, we denote
	\begin{align}
	\underset{k\in \mathbb{N}}{\text{inf}}\alpha _{i,k}=:\alpha _i>0.
	\end{align}
\end{enumerate}

\begin{rem} 
	{\rm For more comments about the the operators $ A_i $ and $ Q_i $ of the assumption (A2), reader can read \cite{cerrai2009khasminskii}.}
\end{rem}

About the coefficients of the system (\ref{orginal1}), we assume it satisfy the following conditions.
\begin{enumerate}[({A}3)]
	\item 
	\begin{enumerate}
		\item The mappings $b_1: \bar{\mathcal{O}}\times \mathbb{R}^2\rightarrow \mathbb{R} $ is continuous and there exists $ m_1\geq 1 $ such that
		\begin{eqnarray}\label{en1}
		\underset{  \xi   \in   \bar{\mathcal{O}}  }{\text{sup}}\left| b_1\left(  \xi ,x,y \right) \right|\le C\left( 1+\left| x \right|^{m_1}+\left| y \right| \right) ,\quad\left( x,y \right) \in \mathbb{R}^2.
		\end{eqnarray}
		\item There exists $C>0$ such that, for any $ x,h\in \mathbb{R}^2 $,
		\begin{eqnarray}\label{en04}
		\underset{  \xi  \in    \bar{\mathcal{O}} }{\text{sup}}\left( b_1\left(  \xi ,x+h \right) -b_1\left( \xi ,x \right) \right) h_1\le C\left| h_1 \right|\left( 1+\left| x \right|+\left| h \right| \right).
		\end{eqnarray}
		\item There exists $ \kappa>0 $ such that
		\begin{eqnarray}\label{en03}
		\underset{  \xi \in  \bar{\mathcal{O}}}{\text{sup}}\left| b_1\left(  \xi ,x \right) -b_1\left( \xi ,y \right) \right|\le C\big( 1+\left| x \right|^{\kappa}+\left| y \right|^{\kappa} \big) \left| x-y \right|,\quad x,y\in \mathbb{R}^2.
		\end{eqnarray}
	\end{enumerate}
\end{enumerate}	
\begin{enumerate}[({A}4)]
	\item 
	\begin{enumerate}
		\item The mappings $b_2:\mathbb{R}\times\bar{\mathcal{O}}\times \mathbb{R}^2\rightarrow \mathbb{R} $ is continuous and there exists $ m_2\geq 1 $ such that
		\begin{eqnarray}\label{en2}
		\underset{\left(t, \xi\right)  \in   \mathbb{R}\times\bar{\mathcal{O}} }{\text{sup}}\left| b_2\left( t,\xi ,x,y \right) \right|\le C\left( 1+\left| x \right|+\left| y \right|^{m_2} \right) ,\quad\left( x,y \right) \in \mathbb{R}^2.
		\end{eqnarray}
		\item There exists $C>0$ such that, for any $ x,h\in \mathbb{R}^2 $,
		\begin{eqnarray}\label{en05}
		\underset{\left(t, \xi\right)  \in   \mathbb{R}\times\bar{\mathcal{O}}}{\text{sup}}\left( b_2\left( t,\xi ,x+h \right) -b_2\left( t,\xi ,x \right) \right) h_2\le C\left| h_2 \right|\left( 1+\left| x \right|+\left| h \right| \right).
		\end{eqnarray}
		\item The mapping $ b_2\left(t,\xi,\cdot \right):\mathbb{R}^2\rightarrow\mathbb{R}  $ is locally Lipschitz-continuous, uniformly with respect to $ \left(t, \xi\right)  \in   \mathbb{R}\times\bar{\mathcal{O}} $. 
		Moreover, for all $\left( t,\xi \right) \in \mathbb{R}\times \bar{\mathcal{O}}$, there exist some continuous function $\varrho:\mathbb{R}\times \bar{\mathcal{O}}\times \mathbb{R}^3\rightarrow \mathbb{R}$
		\begin{eqnarray}\label{en27}
		b_2\left( t,\xi ,x_1,y \right) -b_2\left( t,\xi ,x_2,y \right) =\varrho \left( t,\xi ,x_1,x_2,y \right),
		\end{eqnarray}
		such that 
		\begin{align}
		\underset{\begin{array}{c}
			\left( t,\xi \right) \in \mathbb{R}\times \bar{\mathcal{O}}\\
			\left( x,y \right) \in \mathbb{R}^2,h>0\\
			\end{array}}{\text{inf}}\varrho \left( t,\xi ,x,x+h,y \right) \underset{\begin{array}{c}
			\left( t,\xi \right) \in \mathbb{R}\times \bar{\mathcal{O}}\\
			\left( x,y \right) \in \mathbb{R}^2,h>0\\
			\end{array}}{\text{sup}}\varrho \left( t,\xi ,x,x+h,y \right) \geq 0,
		\end{align}
		and for any $ R>0 $, there exists $ L_R>0, $ such that 
		\begin{eqnarray}\label{216}
		x_1,x_2\in B_{\mathbb{R}}\left( R \right) \Rightarrow \underset{  
			y\in \mathbb{R},\  
		 \left( t,\xi \right) \in \mathbb{R}\times \bar{\mathcal{O}} 
			 }{\sup}\left|\varrho \left( t,\xi ,x_1,x_2,y \right)  \right|\le L_R\left| x_1-x_2 \right|.
		\end{eqnarray}
		\item	
		For all $x,y_1,y_2\in \mathbb{R}$, we have
		\begin{eqnarray}\label{en29}
		\begin{split}
		b_2\left( t,\xi ,x,y_1 \right) -b_2\left( t,\xi ,x,y_2 \right) =-\tau \left( t,\xi ,x,y_1,y_2 \right) \left( y_1-y_2 \right).
		\end{split}
		\end{eqnarray}
		for some measurable function $\tau :\mathbb{R}\times \bar{\mathcal{O}}\times \mathbb{R}^3\rightarrow \left[ 0,\infty \right)$.
	\end{enumerate}
\end{enumerate}
\begin{rem} 
	{\rm Articles \cite{cerrai2011averaging} and \cite{cerrai2017averaging} gives specific examples, which satisfy the conditions of assumptions (A3) and (A4). Readers can read it, and we will not give details in this paper.}
\end{rem} 
\begin{enumerate}[({A}5)]
	\item 
	\begin{enumerate}
		\item The mappings $f_1: \bar{\mathcal{O}}\times \mathbb{R}\rightarrow \mathbb{R}, g_1: \bar{\mathcal{O}}\times \mathbb{R}\times \mathbb{Z}\rightarrow \mathbb{R}, f_2:\mathbb{R}\times\bar{\mathcal{O}}\times \mathbb{R} \rightarrow \mathbb{R}, g_2:\mathbb{R}\times\bar{\mathcal{O}}\times \mathbb{R} \times \mathbb{Z}\rightarrow \mathbb{R}$ are continuous, and the mappings $f_1\left( \xi ,\cdot \right):\mathbb{R}\rightarrow \mathbb{R},g_1\left( \xi,\cdot,z \right):\mathbb{R}\rightarrow \mathbb{R},  f_2\left(t, \xi ,\cdot \right):\mathbb{R} \rightarrow \mathbb{R},g_2\left( t,\xi,\cdot,z \right):\mathbb{R} \rightarrow \mathbb{R}$ are Lipschitz-continuous, uniformly with respect to $(t,\xi,z) \in \mathbb{R}\times \bar{\mathcal{O}}\times \mathbb{Z}$. 
		Moreover, for all  $ p\ge 1 $, there exist positive constants $C_1,C_2 $,  such that for all $ x, y \in \mathbb{R} $, have
		\begin{align}
		\underset{ \xi \in  \bar{\mathcal{O}}}{\text{sup}}\int_{\mathbb{Z}}{\left|  g_1\left(  \xi ,x,z \right) -g_1\left(  \xi ,y,z \right)\right|  }^{p}\upsilon _1\left( dz \right) \leq C_1\left|  x-y\right| ^{p}. 
		\end{align}
		and
		\begin{align}
		\underset{\left( t,\xi \right) \in \mathbb{R}\times \bar{\mathcal{O}}}{\text{sup}}\int_{\mathbb{Z}}{\left|  g_2\left(t, \xi ,x,z \right) -g_2\left(t, \xi ,y,z \right)\right|  }^{p}\upsilon _2\left( dz \right) \leq C_2\left|  x-y\right| ^{p}. 
		\end{align}
		\item For any $  x,y \in \mathbb{R}, $ it hold that
		\begin{eqnarray}
		\underset{ \xi \in  \bar{\mathcal{O}}}{\text{sup}}\Big(\left| f_1\left( \xi ,x \right) \right|^p+\int_{\mathbb{Z}}{\left| g_1\left( \xi ,x, z \right) \right|^p} Y^1\left( dz \right) \Big)\le C\big( 1+\left| x \right|^{\frac{p}{m_1}} \big), 
		\end{eqnarray}
		and
		\begin{eqnarray}
		\underset{\left( t,\xi \right) \in \mathbb{R}\times \bar{\mathcal{O}}}{\text{sup}}\Big(\left| f_2\left( t,\xi ,x \right) \right|^p+\int_{\mathbb{Z}}{\left| g_2\left( t,\xi ,x, z \right) \right|^p} Y^2\left( dz \right) \Big)\le C\big( 1+\left| x \right|^{\frac{p}{m_2}} \big), 
		\end{eqnarray}
		where $ m_1 $ and $ m_2 $ are the constants introduced in (\ref{en1}) and (\ref{en2}).
	\end{enumerate}
\end{enumerate}
\begin{rem}\label{rem2.1}
	{\rm For any $\left( t,\xi \right) \in \mathbb{R}\times \bar{\mathcal{O}}$ and $x,y,h\in   E ,z\in \mathbb{Z}$, we shall set
		\begin{gather}
		B_1\left(  x,y \right) \left( \xi \right) :=b_1\left(  \xi ,x\left( \xi \right) ,y\left( \xi \right) \right) , \quad B_2\left( t,x,y \right) \left( \xi \right) :=b_2\left( t,\xi ,x\left( \xi \right) ,y\left( \xi \right) \right) ,\cr
		\left[ F_1\left(  x \right) h \right] \left( \xi \right) :=f_1\left(  \xi ,x\left( \xi \right) \right) h\left( \xi \right), \quad \left[ F_2\left( t,x  \right) h \right] \left( \xi \right) :=f_2\left( t,\xi ,x\left( \xi \right)  \right) h\left( \xi \right) ,\cr
		\left[ G_1\left(  x,z \right) h \right] \left( \xi \right) :=g_1\left(  \xi ,x\left( \xi \right) ,z \right) h\left( \xi \right), \quad
		\left[ G_2\left( t,x,z \right) h \right] \left( \xi \right) :=g_2\left( t,\xi ,x\left( \xi \right),z \right) h\left( \xi \right).\nonumber
		\end{gather}
		Due to the assumption (A3)  and (A4), we know the mappings
		$ B_1 : E \times  E \rightarrow  E $ and $ B_2: \mathbb{R}\times  E \times   E \rightarrow  E $ are well defined and continuous. According to (\ref{en1}) and (\ref{en2}), for any $ x,y\in E $ and $ t\in\mathbb{R} $, we have 
		\begin{eqnarray}\label{en02}
		\left\|  B_1\left( x,y \right) \right\| \le C\left( 1+\left\|  x \right\|  ^{m_1}+\left\|  y \right\|   \right),\quad \left\|  B_2\left( t,x,y \right) \right\| \le C\left( 1+\left\|  x \right\|  +\left\|  y \right\|^{m_2}  \right).
		\end{eqnarray}
		As a consequence of (\ref{en04}) and (\ref{en05}), it is immediate to check that, for any $ x,y,h,k\in E $, any $ t\in\mathbb{R} $, and any $ \delta\in\mathcal{M} _h, $
		\begin{eqnarray}\label{en07}
		\left< B_1\left( x+h,y+k \right) -B_1\left( x,y \right) ,\delta \right>  \le C\left( 1+\left\| x \right\|  +\left\| y \right\|  +\left\| h \right\|  +\left\| k \right\|   \right), 
		\end{eqnarray}
		and 
		\begin{eqnarray}\label{en06}
		\left< B_2\left( t,x+h,y+k \right) -B_2\left( t,x,y \right) ,\delta \right>  \le C\left( 1+\left\| x \right\|  +\left\| y \right\|  +\left\| h \right\|  +\left\| k \right\|   \right). 
		\end{eqnarray}
		In view of (\ref{en03}), for any $x_1,y_1,x_2,y_2\in E $, we have 
		\begin{small}
			\begin{align}\label{en0221}
			 \left\| B_1\left(  x_1,y_1 \right) -B_1\left( x_2,y_2 \right) \right\| \le C\big( 1+\left\| \left( x_1,y_1 \right) \right\| _{ E \times  E }^{\kappa}+\left\| \left( x_2,y_2 \right) \right\| _{ E \times  E }^{\kappa} \big) \left( \left\| x_1-x_2 \right\| +\left\| y_1-y_2 \right\|  \right). 
			\end{align}
		\end{small}
		In addition, from the equation (\ref{en29}), for every $ \delta\in\mathcal{M} _k, $ we have 
		\begin{align}\label{en223}
		\left\langle B_2\left(t,x,y+k \right)-B_2\left( t,x,y \right),\delta \right\rangle  \leq0.
		\end{align}

		Due to  the assumptions (A5) and  (A6), for any fixed $\left( t,z\right) \in\left( \mathbb{R},\mathbb{Z}\right) ,$ the mappings
		\begin{gather}
		F_1\left( \cdot \right) :  E \rightarrow \mathcal{L}\left(  E \right),\quad 
		G_1\left( \cdot ,z \right) :  E \rightarrow \mathcal{L}\left(  E \right),\cr 
		F_2\left( t,\cdot \right) :  E \rightarrow \mathcal{L}\left(  E \right),\quad 
		G_2\left( t,\cdot ,z \right) :  E \rightarrow \mathcal{L}\left(  E \right),  \nonumber 
		\end{gather}
		are Lipschitz-continuous.}
\end{rem}

\begin{enumerate}[({A}6)]
	\item 
	\begin{enumerate}
		\item The functions $ \gamma :\mathbb{R}\rightarrow \left( 0,\infty \right)  $  and $ l :\mathbb{R}\times \mathcal{O}\rightarrow \mathbb{R}^d $  are periodic, with the same period.
		\item The families of functions
		\begin{gather}
		\mathbf{B}_{ R}:=\left\{ b_2\left( \cdot ,\xi ,\sigma \right) :\ \xi \in \mathcal{O},\ \sigma \in B_{\mathbb{R}^2}\left( R \right) \right\},\cr
		\mathbf{F}_R:=\left\{ f_2\left( \cdot ,\xi ,\sigma \right) :\ \xi \in \mathcal{O},\ \sigma \in B_{\mathbb{R} }\left( R \right) \right\},\cr
		\mathbf{G}_R:=\left\{ g_2\left( \cdot ,\xi ,\sigma ,z \right) :\ \xi \in \mathcal{O},\ \sigma \in B_{\mathbb{R} }\left( R \right) ,\ z\in \mathbb{Z} \right\},\nonumber
		\end{gather}
		are uniformly almost periodic for any  $ R>0 $. 
	\end{enumerate}
\end{enumerate}
\begin{rem}\label{rem4.5}
	{\rm Similar with the proof of \cite[Lemma 6.2]{cerrai2017averaging}, we get that under the assumption (A6), for any  $ R>0 $, the families of functions
		\begin{gather}
		\left\{ B_2\left( \cdot ,x,y \right) :\left( x,y \right) \in B_{E\times E}\left( R \right) \right\},  \quad
		\left\{ F_2\left( \cdot , y \right) : y \in B_{E }\left( R \right) \right\}, \quad 
		\left\{ G_2\left( \cdot , y,z \right) :\left( y,z \right) \in B_{E }\left( R \right) \times \mathbb{Z} \right\}  \nonumber
		\end{gather}
		are uniformly almost periodic.}
\end{rem}

According to the above introduced, system (\ref{orginal1}) can be rewritten as:
\begin{eqnarray}\label{orginal2}
\begin{split}
\begin{cases}
dX^{\epsilon}\left( t \right) &=\left[ A_1  X^{\epsilon}\left( t \right) +B_1\left(  X^{\epsilon}\left( t \right) , Y^{\epsilon}\left( t \right) \right) \right] dt+F_1\left(  X^{\epsilon}\left( t \right) \right) dW ^{Q_1}\left( t \right) \\
&\quad +\int_{\mathbb{Z}}{G_1\left(  X^{\epsilon}\left( t \right) ,z \right)}\tilde{N}_1\left( dt,dz \right), \\
d Y^{\epsilon}\left( t \right) &=\frac{1}{\epsilon}\left[ \left( A_2\left( t \right) -\alpha \right) Y^{\epsilon}\left( t \right) +B_2\left( t,X^{\epsilon}\left( t \right) , Y^{\epsilon}\left( t \right) \right) \right] dt \\
&\quad +\frac{1}{\sqrt{\epsilon}}F_2\left( t, Y^{\epsilon}\left( t \right) \right) dW ^{Q_2}\left( t \right)  +\int_{\mathbb{Z}}{G_2}\left( t, Y^{\epsilon}\left( t \right)  ,z\right) \tilde{N}_{2}^{\epsilon}\left( dt,dz \right), \\
X^{\epsilon}\left( 0 \right) &=x, \quad Y^{\epsilon}\left( 0 \right) =y.
\end{cases}
\end{split}
\end{eqnarray}
 
 We define
$$
\gamma  \left( t,s \right) :=\int_s^t{\gamma  \left( r \right)}dr, \quad s<t.
$$
For any $\epsilon >0$ and $\beta \ge 0$, set
$$
U_{\beta ,\epsilon }\left( t,s \right) =e^{\frac{1}{\epsilon}\gamma  \left( t,s \right) A_2-\frac{\beta}{\epsilon}\left( t-s \right)}, \quad s<t,
$$
in the case $\epsilon =1$, we write $U_{\beta }\left( t,s \right)$, and 
in the case $\epsilon =1$ and $\beta =0$, we write $U \left( t,s \right) $.

Next, for any $\epsilon >0,\beta \ge 0$ and for any $X\in \mathbb{D}\big( \left[ s,t \right] ; W^{1,p}_0(\mathcal{O}) \big) , $ we define
$$
\psi _{\beta ,\epsilon }\left( X;s \right) \left( r \right) =\frac{1}{\epsilon}\int_s^r{U_{\beta ,\epsilon }\left( r,\rho \right) L \left( \rho \right) X\left( \rho \right)}d\rho, \quad s<r<t,
$$
in the case $\epsilon =1$, we write $\psi _{\beta }\left( X;s \right) \left( r \right)$, and in the case $\epsilon =1$ and $\beta =0$, we write $\psi  \left( X;s \right) \left( r \right)$. Then, we can easily get that $ \psi _{\beta ,\epsilon }\left( X;s \right) \left( t \right)  $ is the solution of
$$
dX\left( t \right) =\frac{1}{\epsilon}\left( A_2\left( t \right) -\beta \right) X\left( t \right) dt, \quad t>s, \ X\left( s \right) =0.
$$

Now, we present the main result of this paper:
\begin{thm}\label{th5.3} Under the assumptions {\rm (A1)-(A6)}, there exists $ \bar{\theta}>0 $,  for any $p\ge 1 $ and any initial value $ x\in D\big(\left(-A_1 \right)^\theta  \big) $ with $\theta\in [0,\bar{\theta}\land 1/p),   $ we have 
	\begin{align}\label{510} 
	\underset{\epsilon \rightarrow 0}{\lim}\mathbb{E}\Big( \underset{t\in \left[ 0,T  \right]}{{\sup}}\left\|  X^{\epsilon}\left( t \right) -\bar{X}\left( t \right) \right\| ^p\Big) =0,
	\end{align}
	where  $ \bar{X} $ is the solution of the averaged equation (see equation (\ref{en67}) below). 
\end{thm}

\section{Existence, uniqueness of the solutions}\label{sec-3}
In this section, the existence and uniqueness of solutions is proved by constructing an auxiliary equation. 
\begin{defn}\label{defn3.1}
	For any fix $ x $ and $y $, process $ \left( X^\epsilon\left( t\right), Y^\epsilon\left( t\right)\right)   $ is called a mild solution of the equation (\ref{orginal2}), if 
	\begin{eqnarray}\label{eq32}
	\begin{split}
	\begin{cases}
	X^{\epsilon}\left( t \right) &=e^{A_1t}x +\int_0^t{e^{A_1\left( t-r\right) } B_1\left(  X^{\epsilon}\left( r \right) , Y^{\epsilon}\left( r \right) \right)}dr
	+\int_0^t{e^{A_1\left( t-r\right)} F_1\left(  X^{\epsilon}\left( r \right) \right)}dW^{Q_1}\left( r \right) \cr
	&\quad+\int_0^t{\int_{\mathbb{Z}}{e^{A_1\left( t-r\right)}G_1\left(  X^{\epsilon}\left( r \right) ,z \right)}}\tilde{N}_1\left( dr,dz \right), \\
	 Y^{\epsilon}\left( t \right) &=U_{\alpha ,\epsilon }\left( t,0 \right) y+\psi _{\alpha ,\epsilon }\left( Y^{\epsilon};0 \right)\left( t \right) +\frac{1}{\epsilon}\int_0^t{U_{\alpha ,\epsilon }\left( t,r \right) B_2\left( r,X^{\epsilon}\left( r \right) , Y^{\epsilon}\left( r \right) \right)}dr \cr
	&\quad +\frac{1}{\sqrt{\epsilon}}\int_0^t{U_{\alpha ,\epsilon }\left( t,r \right) F_2\left( r, Y^{\epsilon}\left( r \right) \right)}dW^{Q_2}\left( r \right) \cr
	&\quad +\int_0^t{\int_{\mathbb{Z}}{U_{\alpha ,\epsilon }\left( t,r \right) G_2\left( r, Y^{\epsilon}\left( r \right) ,z \right)}}\tilde{N}_{2}^{\epsilon}\left( dr,dz \right).
	\end{cases} 
	\end{split}
	\end{eqnarray}
\end{defn}

For any $ n\in \mathbb{N} $ and $ \sigma \in \mathbb{R}^2 $, we define 
\begin{eqnarray}
b_{1,n}\left( \xi ,\sigma \right) :=\left\{ \begin{array}{c}
b_1\left( \xi ,\sigma \right) ,\\
b_1\left(  \xi ,\sigma n/\left| \sigma \right| \right) ,\\
\end{array} \right. \quad \begin{array}{c}
\text{if} \ \left| \sigma \right|\leq n,\\
\text{if} \ \left| \sigma \right|>n,\\
\end{array}
\end{eqnarray}
and 
\begin{eqnarray}
b_{2,n}\left( t,\xi ,\sigma \right) :=\left\{ \begin{array}{c}
b_2\left( t,\xi ,\sigma \right) ,\\
b_2\left( t,\xi ,\sigma n/\left| \sigma \right| \right) ,\\
\end{array} \right. \quad \begin{array}{c}
\text{if}\ \left| \sigma \right|\leq n,\\
\text{if}\ \left| \sigma \right|>n.\\
\end{array}
\end{eqnarray}
For each $ b_{i,n} $, denote the corresponding composition operator is $ B_{i,n}. $
It is easy to get that the mapping $ B_{1,n}\left(  \cdot \right) $ and $ B_{2,n}\left( t,  \cdot \right) \left(t \in  \mathbb{R} \right)  $ are Lipschitz-continuous and the mapping $ B_{1,n} $ and $ B_{2,n} $ satisfy all conditions in (A3) and (A4) respectively.  Moreover, for any $  \left( x,y\right) \in E \times E, $ if $ m\leq n, $ we can know that
\begin{align}\label{en11}
\left\| \left( x,y\right) \right\| _{E\times E}\leq m\Rightarrow B_{1,m}\left(  x,y \right) =B_{1,n}\left(  x,y \right) =B_1\left( x,y\right), \quad t\in \mathbb{R},
\end{align}
\begin{align}\label{en10}
\left\| \left( x,y\right) \right\| _{E\times E}\leq m\Rightarrow B_{2,m}\left( t,x,y \right) =B_{2,n}\left( t,x,y \right) =B_2\left( t,x,y\right), \quad t\in \mathbb{R}. 
\end{align}

In order to prove the existence and uniqueness of solutions for systems (\ref{orginal2}), we construct and research the following equation firstly:
\begin{eqnarray}\label{orginal3}
\begin{split}
\begin{cases}
dX^{\epsilon}\left( t \right) &=\left[ A_1  X^{\epsilon}\left( t \right) +B_{1,n}\left(  X^{\epsilon}\left( t \right) , Y^{\epsilon}\left( t \right) \right) \right] dt+F_1\left(  X^{\epsilon}\left( t \right) \right) dW ^{Q_1}\left( t \right) \\
&\quad +\int_{\mathbb{Z}}{G_1\left(  X^{\epsilon}\left( t \right) ,z \right)}\tilde{N}_1\left( dt,dz \right), \\
d Y^{\epsilon}\left( t \right) &=\frac{1}{\epsilon }\left[ \left( A_2\left( t \right) -\alpha \right) Y^{\epsilon}\left( t \right) +B_{2,n}\left( t,X^{\epsilon}\left( t \right) , Y^{\epsilon}\left( t \right) \right) \right] dt \\
&\quad +\frac{1}{\sqrt{\epsilon }}F_2\left( t, Y^{\epsilon}\left( t \right) \right) dW ^{Q_2}\left( t \right)  +\int_{\mathbb{Z}}{G_2}\left( t, Y^{\epsilon}\left( t \right)  ,z\right) \tilde{N}_{2}^{\epsilon }\left( dt,dz \right), \\
X^{\epsilon}\left( 0 \right) &=x, \quad Y^{\epsilon}\left( 0 \right) =y.
\end{cases}
\end{split}
\end{eqnarray}
Due to the coefficients $ B_{i,n}, F_i, G_i $ are Lipschitz continuous, it is easy to prove that the equation (\ref{orginal3}) has solution $ \left(X^{\epsilon,n}\left( t \right), Y^{\epsilon,n}\left( t \right)\right) $ \cite{cerrai2003stochastic,peszat2007stochastic,pei2017two}. 

Next, by proceeding as Lemma 3.1 in \cite{cerrai2011averaging} and Lemma 3.1 in \cite{Xu2018Averaging}, we can prove that the solution $ X^{\epsilon,n}\left( t \right)  $ and $   Y^{\epsilon,n}\left( t \right) $ of equation (\ref{orginal3}) are bounded, the detailed proof will be provided in the Appendix.
\begin{lem}\label{lem3.3}
	Under the assumptions {\rm (A1)-(A5)}, for any $p\ge 1$, there exists a positive constant $C_{p,T}$, such that for any $\epsilon \in \left( 0,1 \right]$, we have   
	\begin{eqnarray}\label{en33}
	\mathbb{E}\underset{t\in \left[ 0,T \right]}{\sup}\left\| X^{\epsilon,n}\left( t \right) \right\|^{p}\le C_{p,T}\left( 1+\left\| x \right\|^{p}+\left\| y \right\|^{p} \right),
	\end{eqnarray}
	and
	\begin{eqnarray}\label{en34}
	\int_0^T{\mathbb{E}\lVert  Y^{ \epsilon,n}\left( r \right) \rVert ^p}dr\le C_{p,T}\left( 1+\left\| x \right\|^{p}+\left\| y \right\|^{p} \right).
	\end{eqnarray}
\end{lem}

Now, we prove the existence and uniqueness of solutions for systems (\ref{orginal2}) through the sequence $ \left\lbrace X^{\epsilon,n} \right\rbrace  $ and $ \left\lbrace Y^{\epsilon,n} \right\rbrace.  $

\begin{thm}\label{th3.2}
	Under the assumptions {\rm (A1)-(A5)}, for any initial value $ x,y $ and $p\ge 1$, there exists unique mild solution $  X^\epsilon $ and $ Y^\epsilon  $ in $ L^p\left( \Omega ;\mathbb{D}\left( \left[ 0,T \right] ;  E  \right)   \right) $ for system (\ref{orginal2}).  Moreover, there exists a positive constant $C_{p,T}$, such that for any $\epsilon \in \left( 0,1 \right]$, we have   
	\begin{eqnarray}\label{en31}
	\mathbb{E}\underset{t\in \left[ 0,T \right]}{\sup}\left\| X^{\epsilon}\left( t \right) \right\|^{p}\le C_{p,T}，
	\left( 1+\left\| x \right\|^{p}+\left\| y \right\|^{p} \right),
	\end{eqnarray}
and
	\begin{eqnarray}\label{en32}
	\int_0^T{\mathbb{E}\lVert  Y^{ \epsilon }\left( r \right) \rVert ^p}dr\le C_{p,T} 
	\left( 1+\left\| x \right\|^{p}+\left\| y \right\|^{p} \right).
	\end{eqnarray}
\end{thm}
\para{Proof:} 
 Fixed $  \epsilon\in (0,1], $ use the same arguement as the equation (\ref{en33}) for $ Y^{\epsilon, n}\left( t \right), $  it is easy to prove that the following equation also holds 
\begin{eqnarray}\label{323} 
\mathbb{E}\underset{t\in \left[ 0,T \right]}{\sup}\left\| Y^{\epsilon, n}\left( t \right) \right\| ^{p}\le C_{p,\epsilon,T} \left( 1+\left\| x \right\| ^{p}+\left\| y \right\| ^{p} \right). 
\end{eqnarray} 
To prove  Theorem \ref{th3.2}, 
for any $ n\in \mathbb{N} $ and fixed $  \epsilon\in (0,1], $  we define
$$
\tau _{n}^{x,y} :=\text{inf}\left\{ t\geq 0:  \left\| X^{\epsilon, n}\left( t \right) \right\|  +  \lVert  Y^{\epsilon, n}\left( t \right) \rVert    \ge n  \right\},  
$$
and let
$$
\tau^{x,y}   :=\underset{n\in \mathbb{N}}{\text{sup}}\tau _{n}^{x,y} .
$$
Due to (\ref{en33}) and (\ref{323}), we can know that the sequence of stopping times $ \left\lbrace \tau_{n}^{x,y} \right\rbrace  $ is non-decreasing and  
 $ \mathbb{P}(\tau^{x,y}   =+\infty)=1 $. Indeed,
 $$
 \mathbb{P}\left( \tau ^{x,y}<+\infty \right) =\lim_{T\rightarrow +\infty}\mathbb{P}\left( \tau ^{x,y}\leq T \right), 
 $$
and for each $ T > 0, $ we have
\begin{align}
 \mathbb{P}\left( \tau ^{x,y}\leq T \right) &=\lim_{n\rightarrow +\infty}\mathbb{P}\left( \tau _{n}^{x,y}\leq T \right)=\lim_{n\rightarrow +\infty}\mathbb{P}\Big( \sup_{t\in \left[ 0,T \right]}\lVert X^{\epsilon, n} \left( t \right) \rVert +\sup_{t\in \left[ 0,T \right]}\lVert Y^{\epsilon, n} \left( t \right) \rVert \geq n \Big)\cr
 &=\lim_{n\rightarrow +\infty}\frac{1}{n^2}\mathbb{E}\Big( \sup_{t\in \left[ 0,T \right]}\lVert X^{\epsilon, n} \left( t \right) \rVert ^2+\sup_{t\in \left[ 0,T \right]}\lVert Y^{\epsilon, n} \left( t \right) \rVert ^2 \Big) \cr
 &=\lim_{n\rightarrow +\infty}C_{p,\epsilon,T}\frac{1+\lVert x \rVert ^2+\lVert y \rVert ^2}{n^2}=0. \nonumber
\end{align}
Hence, $ P(\tau^{x,y}   =+\infty)=1. $ Then, for any $ t\in [0,T] $ and $ \omega\in \left\lbrace \tau^{x,y}   =+\infty\right\rbrace  $, there exists $ m\in \mathbb{N} $ such that $ t\leq \tau _{m}^{x,y} (\omega)$. Then, we define
\begin{eqnarray}\label{en8}
 X^{\epsilon}  (t)(\omega):= X^{\epsilon, m}  (t)(\omega) \quad \text{and} \quad
 Y^{\epsilon}   (t)(\omega):= Y^{\epsilon, m}  (t)(\omega).
\end{eqnarray}
For any $ n\geq m $, set $ \tau:=\tau_{n}^{x,y}\land\tau_{m}^{x,y}. $  Thanks to (\ref{en11}) and (\ref{en10}), for any $ t\leq\tau, $ we can get
\begin{align}\label{enn2}
&\qquad X^{\epsilon, n}\left( t\land \tau \right) - X^{\epsilon, m}\left( t\land \tau \right)+ Y^{\epsilon, n}\left( t\land \tau \right) - Y^{\epsilon, m}\left( t\land \tau \right) \cr
&=\int_0^{t\land \tau}{e^{A_1\left( t\land \tau -r \right)}\left[ B_{1,n}\left( X^{\epsilon, n}\left( r \right) , Y^{\epsilon, n}\left( r \right) \right) -B_{1,m}\left( X^{\epsilon, m}\left( r \right) , Y^{\epsilon, m}\left( r \right) \right) \right]}dr\cr
&\quad+\int_0^{t\land \tau}{e^{A_1\left( t\land \tau -r \right)}\left[ F_1\left( X^{\epsilon, n}\left( r \right) \right) -F_1\left( X^{\epsilon, m}\left( r \right) \right) \right]}dW^{Q_1}\left( r \right) \cr
&\quad+\int_0^{t\land \tau}{\int_{\mathbb{Z}}{e^{A_1\left( t\land \tau -r \right)}\left[ G_1\left( X^{ \epsilon,n}\left( r \right) ,z \right) -G_1\left( X^{\epsilon, m}\left( r \right) ,z \right) \right]}}\tilde{N}_1\left( dr,dz \right) \cr
&\quad+\psi _{\alpha,\epsilon }\left( Y^{\epsilon, n};0 \right) \left( t\land \tau \right) -\psi _{\alpha,\epsilon }\left( Y^{\epsilon, m};0 \right) \left( t\land \tau \right) \cr
&\quad+\frac{1}{\epsilon}\int_0^{t\land \tau}{U_{\alpha ,\epsilon }\left( t\land\tau,r \right)\left[ B_{2,n}\left( r, X^{\epsilon, n}\left( r \right) , Y^{\epsilon, n}\left( r \right) \right) -B_{2,m}\left( r, X^{\epsilon, m}\left( r \right) , Y^{\epsilon, m}\left( r \right) \right) \right]}dr\cr
&\quad+\frac{1}{\sqrt{\epsilon}}\int_0^{t\land \tau}{U_{\alpha ,\epsilon }\left( t\land\tau,r \right)\left[ F_2\left( r, Y^{ \epsilon,n}\left( r \right) \right) -F_2\left( r, Y^{\epsilon, m}\left( r \right) \right) \right]}dW^{Q_2}\left( r \right) \cr
&\quad+\int_0^{t\land \tau}{\int_{\mathbb{Z}}{U_{\alpha ,\epsilon }\left( t\land\tau,r \right)\left[ G_2\left( r, Y^{ \epsilon,n}\left( r \right) ,z \right) -G_2\left( r, Y^{\epsilon, m}\left( r \right) ,z \right) \right]}}\tilde{N}_{2}^\epsilon \left( dr,dz \right) \cr
&=\int_0^t{I_{\left\{ r\le \tau \right\}}e^{A_1\left( t\land \tau -r \right)}\left[ B_{1,n}\left( X^{\epsilon, n}\left( r\land \tau \right) , Y^{\epsilon, n}\left( r\land \tau \right) \right) -B_{1,n}\left( X^{\epsilon, m}\left( r\land \tau \right) , Y^{\epsilon, m}\left( r\land \tau \right) \right) \right]}dr\cr
&\quad+\int_0^t{I_{\left\{ r\le \tau \right\}}e^{A_1\left( t\land \tau -r \right)}\left[ F_1\left( X^{\epsilon, n}\left( r\land \tau \right) \right) -F_1\left( X^{\epsilon, m}\left( r\land \tau \right) \right) \right]}dW^{Q_1}\left( r \right) \cr
&\quad+\int_0^t{\int_{\mathbb{Z}}{I_{\left\{ r\le \tau \right\}}e^{A_1\left( t\land \tau -r \right)}\left[ G_1\left( X^{\epsilon, n}\left( r\land \tau \right) ,z \right) -G_1\left( X^{\epsilon, m}\left( r\land \tau \right) ,z \right) \right]}}\tilde{N}_1\left( dr,dz \right) \cr
&\quad+\psi _{\alpha,\epsilon  }\left( Y^{\epsilon,n}- Y^{\epsilon, m};0 \right) \left( t\land \tau \right)+\frac{1}{\epsilon}\int_0^t{I_{\left\{ r\le \tau \right\}}U_{\alpha ,\epsilon }\left( t\land\tau,r \right)\big[ B_{2,n}\left( r\land \tau , X^{\epsilon, n}\left( r\land \tau \right) , Y^{\epsilon, n}\left( r\land \tau \right) \right)} \cr
&\qquad\qquad\qquad\qquad\qquad\qquad\qquad\qquad\qquad\qquad\qquad\qquad -B_{2,n}\left( r\land \tau , X^{\epsilon, m}\left( r\land \tau \right) , Y^{\epsilon, m}\left( r\land \tau \right) \right) \big]dr\cr
&\quad+\frac{1}{\sqrt{\epsilon}}\int_0^t{I_{\left\{ r\le \tau \right\}}U_{\alpha ,\epsilon }\left( t\land\tau,r \right)\left[ F_2\left( r\land \tau , Y^{\epsilon, n}\left( r\land \tau \right) \right) -F_2\left( r\land \tau , Y^{\epsilon, m}\left( r\land \tau \right) \right) \right]}dW^{Q_2}\left( r \right) \cr
&\quad+\int_0^t{\int_{\mathbb{Z}}{I_{\left\{ r\le \tau \right\}}U_{\alpha ,\epsilon }\left( t\land\tau,r \right)\left[ G_2\left( r\land \tau , Y^{\epsilon, n}\left( r\land \tau \right) ,z \right) -G_2\left( r\land \tau , Y^{\epsilon, m}\left( r\land \tau \right) ,z \right) \right]}}\tilde{N}_{2}^\epsilon \left( dr,dz \right) \cr
&:=\sum_{i=1}^7{\mathcal{I}_i\left( t \right)}.
\end{align}
For $\mathcal{I}_{2}\left( t\right)+\mathcal{I}_6\left(t \right) ,  $ according to equation (3.5) in \cite{cerrai2009khasminskii},
fixed  $ \bar{p}>1 $, such that  $ \frac{\beta _1\left( \rho _1-2 \right)}{\rho _1}\frac{\bar{p}}{\bar{p}-2}\vee\frac{\beta _2\left( \rho _2-2 \right)}{\rho _2}\frac{\bar{p}}{\bar{p}-2}<1 $. Then, for any  $p\ge \bar{p},$ we can get
\begin{align}
&\qquad\mathbb{E}\lVert \mathcal{I}_2\left( t \right) \lVert ^p+\mathbb{E} \rVert \mathcal{I}_6\left( t \right) \rVert ^p\cr
&\le C_p\Big( \int_0^t{\lVert e^{A_1\left( t\land \tau -r \right)}\left[ F_1\left( X^{\epsilon ,n}\left( r\land \tau \right) \right) -F_1\left( X^{\epsilon ,m}\left( r\land \tau \right) \right) \right] I_{\left\{ r\le \tau \right\}}Q_1 \rVert _{2}^{2}}dr \Big) ^{\frac{p}{2}}\cr
&\quad+C_{p,\epsilon}\Big( \int_0^t{\lVert U_{\alpha ,\epsilon}\left( t\land \tau ,r \right) \left[ F_2\left( r\land \tau ,Y^{\epsilon ,n}\left( r\land \tau \right) \right) -F_2\left( r\land \tau ,Y^{\epsilon ,m}\left( r\land \tau \right) \right) \right] I_{\left\{ r\le \tau \right\}}Q_2 \rVert _{2}^{2}}dr \Big) ^{\frac{p}{2}}\cr
&\le C_p\mathbb{E}\Big( \int_0^t{\left( t\land \tau -r \right) ^{-\frac{\beta _1\left( \rho _1-2 \right)}{\rho _1}}e^{-\frac{\alpha _1\left( \rho _1+2 \right)}{\rho _1}\left( t\land \tau -r \right)}I_{\left\{ r\le \tau \right\}}\lVert X^{\epsilon ,n}\left( r\land \tau \right) - X^{\epsilon ,m}\left( r\land \tau \right) \rVert ^2}dr \Big) ^{\frac{p}{2}}\cr
&\quad+C_{p,\epsilon}\mathbb{E}\Big( \int_0^t{\left( {\gamma \left( t\land\tau,r \right)}/{\epsilon}
	 \right) ^{-\frac{\beta _2\left( \rho _2-2 \right)}{\rho _2}}e^{-\frac{\alpha _2\left( \rho _2+2 \right)}{\rho _2}\left( {\gamma \left( t\land\tau,r \right)}/{\epsilon}
	 	\right)}I_{\left\{ r\le \tau \right\}}\lVert Y^{\epsilon ,n}\left( r\land \tau \right) - Y^{\epsilon ,m}\left( r\land \tau \right) \rVert ^2}dr \Big) ^{\frac{p}{2}}\cr
&\le C_p\mathbb{E}\Big[\Big( \int_0^{t\land \tau}{\left( t\land \tau -r \right) ^{-\frac{\beta _1\left( \rho _1-2 \right)}{\rho _1}\frac{p}{p-2}}}dr \Big) ^{\frac{p-2}{2}}\int_0^t{\lVert X^{\epsilon ,n}\left( r\land \tau \right) - X^{\epsilon ,m}\left( r\land \tau \right) \rVert ^p}dr\Big]\cr
&\quad+C_{p,\epsilon}\mathbb{E}\Big[\Big( \int_0^{t\land \tau}{\left[ \gamma_0\left( t\land \tau -r \right)\right]  ^{-\frac{\beta _2\left( \rho _2-2 \right)}{\rho _2}\frac{p}{p-2}}}dr \Big) ^{\frac{p-2}{2}}\int_0^t{\lVert Y^{\epsilon ,n}\left( r\land \tau \right) - Y^{\epsilon ,m}\left( r\land \tau \right) \rVert ^p}dr\Big] \cr
&\le C_{p,\epsilon,T}\int_0^t{\mathbb{E}\lVert X^{\epsilon ,n}\left( r\land \tau \right) - X^{\epsilon ,m}\left( r\land \tau \right) \rVert ^p+\mathbb{E}\lVert Y^{\epsilon ,n}\left( r\land \tau \right) - Y^{\epsilon ,m}\left( r\land \tau \right) \rVert ^p}dr.
\end{align}
Because $ B_{1,n}\left( \cdot,\cdot\right)   $ and $ B_{2,n}\left( r,\cdot,\cdot\right)  $ are Lipschitz-continuous uniformly with respect to $ t \in \mathbb{R}. $ Then, for any $ t\in [0,T], $ using the H\"{o}lder inequality and Kunita's first inequality \cite[Theorem 4.4.23]{Applebaum2009Processes},  it is easy to get that
\begin{align}
&\qquad\mathbb{E}\left\| \mathcal{I}_1\left( t \right) \right\| ^p+\mathbb{E}\left\| \mathcal{I}_3\left( t \right) \right\| ^p+ \mathbb{E}\left\| \mathcal{I}_5\left( t \right)\right\|^p+\mathbb{E}\left\| \mathcal{I}_7\left( t \right) \right\| ^p  \cr
&\le C_{p,\epsilon,T}\int_0^t{\mathbb{E} \lVert X^{\epsilon, n}\left( r \land \tau \right) - X^{ \epsilon,m}\left( r \land \tau \right) \rVert^p +\mathbb{E} \lVert Y^{ \epsilon,n}\left( r \land \tau \right) - Y^{\epsilon, m}\left( r \land \tau \right) \rVert^p}dr.
\end{align}
Moreover, as a consequence of (2.15) in \cite{cerrai2017averaging}, we know that $ \psi_{\alpha,\epsilon   }\left( \cdot;0 \right)   $ is a bounded linear operator, and there exists a continuous increasing function $ C_\alpha $ with $ C_\alpha \left( 0\right) =0  $ and     $ C_\alpha \left( \cdot\right) \rightarrow 0\left(\text{if} \ \alpha  \ \text{large enough} \right)   $ such that, for any $ t>0, $ have
\begin{align}\label{56}
\mathbb{E}\sup_{  t \in [0,T] } \left\| \mathcal{I}_4 \left( t \right)\right\|  \le    C_\alpha\left( T  \right)  \mathbb{E} \sup_{  t \in [0,T] }\left\| Y^{ \epsilon,n}\left(t \land \tau\right) - Y^{\epsilon, m}\left( t\land \tau\right) \right\|. 
\end{align}
 Hence, for any $ t\in [0,T], $ as $ \alpha  $ is large enough and thanks to (\ref{enn2})-(\ref{56}), using the Gronwall inequality, we can get
\begin{align}\label{en023}
\mathbb{E}\left\| X^{\epsilon, n}\left( t\land \tau \right) - X^{\epsilon, m}\left( t\land \tau \right)\right\|^p +\mathbb{E}\left\| Y^{ \epsilon,n}\left( t\land \tau \right) - Y^{ \epsilon,m}\left( t\land \tau \right) \right\| ^p=0,
\end{align}
it follows that
\begin{align}\label{en025}
 X^{\epsilon ,n}\left( t  \right) = X^{\epsilon ,m}\left( t  \right)   \quad\text{and} \quad  
 Y^{\epsilon ,n}\left( t  \right) = Y^{\epsilon ,m}\left( t  \right), \quad t\leq \tau_{m}^{x,y}\land\tau_{n}^{x,y}.
\end{align}

Next,  thanks to (\ref{en8}),  (\ref{en11}) and (\ref{en10}), if $ \omega\in \left\lbrace \tau   =+\infty\right\rbrace  $ and $ t \leq \tau_m^{x,y}, $ it yields
\begin{eqnarray} 
\begin{split}
\begin{cases}
 X^{\epsilon}\left( t \right) &=e^{A_1t}x +\int_0^t{e^{A_1\left( t-r\right) } B_1\left(  X^{\epsilon}\left( r \right) , Y^{\epsilon}\left( r \right) \right)}dr
+\int_0^t{e^{A_1\left( t-r\right)} F_1\left(  X^{\epsilon}\left( r \right) \right)}dW^{Q_1}\left( r \right) \cr
&\quad+\int_0^t{\int_{\mathbb{Z}}{e^{A_1\left( t-r\right)}G_1\left(  X^{\epsilon}\left( r \right) ,z \right)}}\tilde{N}_1\left( dr,dz \right), \\
 Y^{\epsilon}\left( t \right) &=U_{\alpha ,\epsilon }\left( t,0 \right) y+\psi _{\alpha ,\epsilon }\left( Y^{\epsilon};0 \right)\left( t \right) +\frac{1}{\epsilon}\int_0^t{U_{\alpha ,\epsilon }\left( t,r \right) B_2\left( r, X^{\epsilon}\left( r \right) , Y^{\epsilon}\left( r \right) \right)}dr \cr
&\quad +\frac{1}{\sqrt{\epsilon}}\int_0^t{U_{\alpha ,\epsilon }\left( t,r \right) F_2\left( r, Y^{\epsilon}\left( r \right) \right)}dW^{Q_2}\left( r \right) \cr
&\quad +\int_0^t{\int_{\mathbb{Z}}{U_{\alpha ,\epsilon }\left( t,r \right) G_2\left( r, Y^{\epsilon}\left( r \right) ,z \right)}}\tilde{N}_{2}^{\epsilon}\left( dr,dz \right),
\end{cases} 
\end{split}
\end{eqnarray}
$ \mathbb{P}-a.s. $, that is, $\left(  X^\epsilon \left( t\right), Y^\epsilon \left( t\right)\right)   $ is the mild solution of the system (\ref{orginal2}).

 Denote another solution of system (\ref{orginal2}) is $\left(  X^{\epsilon,*} \left( t\right), Y^{\epsilon,*} \left( t\right)\right),    $ by proceeding as the equation (\ref{en025}), we can also get that 
\begin{eqnarray}
X^{\epsilon,*} \left( t\right) = X^\epsilon  \left( t\right) \quad
\text{and}\quad  Y^{\epsilon,*} \left( t\right) = Y^\epsilon \left( t\right).\nonumber
\end{eqnarray}
Thus, we prove the solution of system (\ref{orginal2}) is unique. 

Finally, using the same argument as \lemref{lem3.3}, we can prove that (\ref{en31}) and (\ref{en32}) hold. 
This completes the proof of \thmref{th3.2}. \qed

\begin{lem}\label{lem3.4}
	Under the assumptions {\rm (A1)-(A5)}, there exists $ \bar{\theta}>0 $,  for any $p\ge 1,\ h \in \left( 0,1\right) $ and any initial value $ x\in D\big(\left(-A_1 \right)^\theta  \big) $ with $\theta\in [0,\bar{\theta}\land 1/p),   $ there exists a positive constant $C_{p, \theta, T}$, such that for any $\epsilon \in \left( 0,1 \right]$, we have   
	\begin{align} 
	\mathbb{E}\left\| X^{\epsilon}\left( t+h \right) - X^{\epsilon}\left( t \right) \right\|^{p}\leq C_{p, \theta, T}h^{\frac{p}{2}\land\theta p\land 1  }\big( 1+\left\| x \right\|^{  m_1 p }+\left\| y \right\|^{m_1 p} \big). 
	\end{align}
\end{lem}
\para{Proof:} From (\ref{eq32}), we have
\begin{align}\label{en318}
&\quad \mathbb{E}\left\| X^{\epsilon}\left( t+h \right) - X^{\epsilon}\left( t \right)\right\| ^p\cr &\leq C_p \big\|\big(e^{A_1 h}-I \big)e^{A_1 t} x\big\| ^p \cr
&\quad+ C_p\mathbb{E}\Big\|\int_0^{t+h}{e^{A_1 (t+h-r)} B_1\left(  X^{\epsilon}\left( r \right) , Y^{\epsilon}\left( r \right) \right) }dr-\int_0^{t }{e^{A_1 (t -r)} B_1\left(  X^{\epsilon}\left( r \right) , Y^{\epsilon}\left( r \right) \right) }dr\Big\| ^p\cr
&\quad+ C_p\mathbb{E}\Big\|\int_0^{t+h}{e^{A_1 (t+h-r)}}F_1\left(  X^{\epsilon}\left( r \right) \right)dW^{Q_1}\left( r \right)-\int_0^{t }{e^{A_1 (t -r)}}F_1\left(  X^{\epsilon}\left( r \right) \right)dW^{Q_1}\left( r \right)\Big\| ^p\cr
&\quad+ C_p\mathbb{E}\Big\|\int_0^{t+h}{\int_{\mathbb{Z}}{e^{A_1 (t+h-r)}}G_1\left(  X^{\epsilon}\left( r \right) ,z \right)}\tilde{N}_1\left( dr,dz \right)-\int_0^{t }{\int_{\mathbb{Z}}{e^{A_1 (t -r)}}G_1\left(  X^{\epsilon}\left( r \right) ,z \right)}\tilde{N}_1\left( dr,dz \right)\Big\| ^p\cr
&:= C_p\sum_{i=1}^4{\mathcal{I}_{t}^{i}}. 
\end{align}

For $ \mathcal{I}_{t}^{1}, $ using the spectral properties \cite{Pazy2012Semigroups} of operator $ A_1, $ for any $ 0<\theta\leq 1, $ there exist some $ C_\theta>0 $ and $ \gamma>0, $ such that
$$
\big\| (-A_1)^\theta e^{A_1 t}\big\| \leq C_\theta t^{-\theta} e^{-\gamma t}, \quad t>0.
$$
Moreover, there exists some positive constant $ C_\theta, $ such that
$$
\big\| (-A_1)^{-\theta} \left( e^{A_1 t}-I\right) \big\| \leq C_\theta t^{ \theta}, \quad t>0.
$$
Therefore, we can get 
\begin{align}\label{en18}
\mathcal{I}_{t}^{1}\leq \big\|\big(e^{A_1 h}-I \big)(-A_1)^{-\theta}\big\| ^p \big\| e^{A_1 t}(-A_1)^{ \theta}\big\| ^p \left\| x\right\|^p \leq C_{p, \theta, T} h^{\theta p}\left\| x\right\|^p
\end{align}

The estimations of $ \mathcal{I}_{t}^{2}  $ is similar and simpler than the estimations of $ \mathcal{I}_{t}^{3} $ and $ \mathcal{I}_{t}^{4}, $ so we only give the following conclusion  
\begin{align}\label{en12}
\mathcal{I}_{t}^{2} \le C_{p,\theta,T}  \big( \Gamma\left( 1-p\theta\right) h^{\theta p}+ h^p \big)\big( 1+\left\| x \right\|^{  m_1 p }+\left\| y \right\|^{m_1 p} \big).
\end{align}

For $ \mathcal{I}_{t}^{3}. $ Thanks to the Lemma 4.1 in \cite{cerrai2009khasminskii}, we can get that
\begin{align}\label{en222}
\lVert e^{A_1\left( t-r \right)}\left( -A_1 \right) ^{\theta}F_1\left( X^{\epsilon}\left( r \right) \right) Q_1 \rVert _{2}^{2}\leq C_{\theta}\left( t-r \right) ^{-\frac{\beta _1\left( \rho _1-2 \right) +\theta \left( \rho _1+2 \right)}{\rho _1}\frac{p}{p-2}}.
\end{align}
Then, if we choose some $ \bar{\theta}>0, $ such that
$$
 \left[ {\beta _1\left( \rho _1-2 \right) +\bar{\theta}\left( \rho _1+2 \right)}\right] /{\rho _1} <1,
$$
and fixed $ \bar{p}>1, $ such that
$$
\frac{\beta _1\left( \rho _1-2 \right) +\bar{\theta}\left( \rho _1+2 \right)}{\rho _1}\frac{\bar{p}}{\bar{p}-2}  <1.
$$
Due to (\ref{en222}) and equation (3.5) in \cite{cerrai2009khasminskii}, for any $ p\geq\bar{p}  $ and $ \theta\in [0,\bar{\theta}], $ we have
\begin{align}\label{en2222}
\mathcal{I}_{t}^{3}&\le C_p\mathbb{E}\Big\|\big( e^{A_1h}-I \big) \left( -A_1 \right) ^{-\theta} \int_0^t{e^{A_1\left( t-r \right)}\left( -A_1 \right) ^{\theta}}F_1\left( X^{\epsilon}\left( r \right) \right) dW^{Q_1}\left( r \right) \Big\| ^p\cr
&\quad+C_p\mathbb{E}\Big\| \int_t^{t+h}{e^{A_1\left( t+h-r \right)}}F_1\left( X^{\epsilon}\left( r \right) \right) dW^{Q_1}\left( r \right) \Big\| ^p\cr
&\le C_p\mathbb{E}\Big( \int_0^t{\lVert e^{A_1\left( t-r \right)}\left( -A_1 \right) ^{\theta}F_1\left( X^{\epsilon}\left( r \right) \right) Q_1 \rVert _{2}^{2}}dr \Big) ^{\frac{p}{2}}\lVert \big( e^{A_1h}-I \big) \left( -A_1 \right) ^{-\theta} \rVert ^p\cr
&\quad+C_p\mathbb{E}\Big( \int_t^{t+h}{\lVert e^{A_1\left( t+h-r \right)}F_1\left( X^{\epsilon}\left( r \right) \right) Q_1 \rVert _{2}^{2}}dr \Big) ^{\frac{p}{2}}\cr
&\le C_{p,\theta}h^{\theta p}\mathbb{E}\Big( \int_0^t{\left( t-r \right) ^{\frac{\beta _1\left( \rho _1-2 \right) +\theta \left( \rho _1+2 \right)}{\rho _1}}\lVert F_1\left( X^{\epsilon}\left( r \right) \right) \rVert ^2}dr \Big) ^{\frac{p}{2}}\cr
&\quad+C_p\mathbb{E}\Big( \int_t^{t+h}{\left( t+h-r \right) ^{-\frac{\beta _1\left( \rho _1-2 \right)}{\rho _1}}e^{-\frac{\alpha _1\left( \rho _1+2 \right)}{\rho _1}\left( t+h-r \right)}\lVert F_1\left( X^{\epsilon}\left( r \right) \right) \rVert ^2}dr \Big) ^{\frac{p}{2}}\cr
&\le C_{p,\theta}h^{\theta p}\Big( \int_0^t{\left( t-r \right) ^{-\frac{\beta _1\left( \rho _1-2 \right) +\theta \left( \rho _1+2 \right)}{\rho _1}\frac{p}{p-2}}}dr \Big) ^{\frac{p-2}{2}}\int_0^t{\mathbb{E}\lVert F_1\left( X^{\epsilon}\left( r \right) \right) \rVert ^p}dr\cr
&\quad+C_p\Big( \int_t^{t+h}{\left( t+h-r \right) ^{-\frac{\beta _1\left( \rho _1-2 \right)}{\rho _1}\frac{p}{p-2}}}dr \Big) ^{\frac{p-2}{2}}\int_t^{t+h}{\mathbb{E}\lVert F_1\left( X^{\epsilon}\left( r \right) \right) \rVert ^p}dr\cr
&\le C_{p,\theta,T}\big( h^{\theta p}+h^{2-\frac{\beta _1\left( \rho _1-2 \right)}{\rho _1}\frac{p}{p-2}} \big) \big( 1+\lVert x \rVert ^{\frac{p}{m_1}}+\lVert y \rVert ^{\frac{p}{m_1}} \big) 
\end{align}

For $ \mathcal{I}_{t}^{4}, $ using Kunita's first inequality and the  H\"{o}lder inequality, thanks to (\ref{en31}), for any $ 0<\theta<1/p, $ we can get
\begin{align}\label{en14}
\mathcal{I}_{t}^{4}&\le C_p\mathbb{E}\Big\|  \int_0^t{\int_{\mathbb{Z}}{\big( e^{A_1 h} -I \big) e^{A_1 (t -r)} G_1\left(  X^{\epsilon}\left( r \right) ,z \right)}\tilde{N}_1\left( dr,dz \right)} \Big\|  ^p\cr
&\quad+C_p\mathbb{E}\Big\| \int_t^{t+h}{\int_{\mathbb{Z}}{e^{A_1 (t+h-r)} G_1\left(  X^{\epsilon}\left( r \right) ,z \right)}\tilde{N}_1\left( dr,dz \right)} \Big\| ^p\cr
&\le C_p\mathbb{E}\Big( \int_0^t{\int_{\mathbb{Z}}{\lVert \big( e^{A_1 h} -I \big) \left( -A_1 \right) ^{-\theta} \rVert ^2\lVert e^{A_1 (t -r)} \left( -A_1 \right) ^{\theta} \rVert ^2\lVert G_1\left(  X^{\epsilon}\left( r \right) ,z \right) \rVert ^2} v_1\left( dz \right) dr} \Big) ^{\frac{p}{2}}\cr
&\quad+C_p\mathbb{E}\int_0^t{\int_{\mathbb{Z}}{\lVert \big( e^{A_1 h} -I \big) \left( -A_1 \right) ^{-\theta} \rVert ^p\lVert e^{A_1 (t -r)} \left( -A_1 \right) ^{\theta} \rVert ^p\lVert G_1\left(   X^{\epsilon}\left( r \right) ,z \right) \rVert ^p} v_1\left( dz \right) dr}\cr
&\quad +C_p\mathbb{E}\Big( \int_t^{t+h}{\int_{\mathbb{Z}}{\lVert G_1\left(  X^{\epsilon}\left( r \right) ,z \right) \rVert ^2} v_1\left( dz \right) dr} \Big) ^{\frac{p}{2}}+C_p\mathbb{E}\int_t^{t+h}{\int_{\mathbb{Z}}{\lVert G_1\left(  X^{\epsilon}\left( r \right) ,z \right) \rVert ^p} v_1\left( dz \right) dr}\cr
&\le C_{p,\theta}\mathbb{E}\Big( \int_0^t{h^{2\theta}\left( t-r\right)^{-2\theta}e^{-2\gamma\left( t-r\right) }\big( 1+\lVert X^{\epsilon}\left( r \right) \rVert ^\frac{2}{m_1} \big) dr} \Big) ^{\frac{p}{2}}\cr
&\quad+C_{p,\theta}\int_0^t{h^{\theta p}\left( t-r\right)^{-p\theta}e^{-p\gamma\left( t-r\right) }\big( 1+\mathbb{E}\lVert X^{\epsilon}\left( r \right) \rVert ^\frac{p}{m_1} \big)}dr\cr
&\quad+C_{p }\mathbb{E}\Big( \int_t^{t+h}{\big( 1+\lVert X^{\epsilon}\left( r \right) \rVert ^{\frac{2}{m_1}} \big) dr} \Big) ^{\frac{p}{2}}+C_p\int_t^{t+h}{\big( 1+\mathbb{E}\lVert X^{\epsilon}\left( r \right) \rVert ^\frac{p}{m_1} \big) dr}\cr
&\le C_{p,\theta,T}  \big( \Gamma\left( 1-p\theta\right) h^{\theta p}+ h^{\frac{p}{2}}+ h \big)\big( 1+\left\| x \right\|^{\frac{p}{m_1}}+\left\| y \right\|^{\frac{p}{m_1}} \big). 
\end{align}
Finally, substituting (\ref{en18})--(\ref{en14}) into (\ref{en318}), for any $p\ge \bar{p} $ and $\theta\in [0,\bar{\theta}\land 1/p),   $ it yields
\begin{align} 
\mathbb{E}\left\| X^{\epsilon}\left( t+h \right) - X^{\epsilon}\left( t \right) \right\|^{p}\leq C_{p, \theta, T}h^{\frac{p}{2}\land\theta p\land 1 }\big( 1+\left\| x \right\|^{  m_1 p }+\left\| y \right\|^{m_1 p} \big), 
\end{align}
Then, using the H\"{o}lder inequality for $p<\bar{p}$,  the proof of \lemref{lem3.4} is complete. \qed

\section{Proof of the main result}\label{sec-4}
In this section, we prove the main \thmref{th5.3}.  
The main research process of this theorem are as follows. Firstly, the averaged equation is defined by studying the existence of an evolution family of measures for the frozen fast equation. Secondly, inspired by Khasminskii in \cite{khas1968on},  using the discretization techniques, we construct an auxiliary process $ \hat{Y}^\epsilon \left( t\right)   $ and derive uniform bounds for it. Thirdly, 
we obtain appropriate control of $ X^{\epsilon}\left( t\right) -\bar{X}\left( t \right)  $ after the stopping time on the basis of some a-priori estimates for the solution $ \left(  X^{\epsilon}\left( t\right), Y^{\epsilon}\left( t\right)\right)  $ of original equation (\ref{orginal1}). Finally, based on the ergodic property of the frozen fast equation, the appropriate control of $ X^{\epsilon}\left( t\right) -\bar{X}\left( t \right)  $ before the stopping time is obtained and the main \thmref{th5.3} is proved. 

\subsection{The averaged equation}
In this part, we main research the fast equation with frozen slow component $ x $. Our target is to prove that there exists an evolution family of measures for this fast equation and define the averaged equation through it.

First, for any $ s\in \mathbb{R} $ and any frozen slow component $x,$ we introduce the following problem
\begin{align}\label{en41}
dY\left( t \right) &=\left[ \left( A_2\left( t \right) -\alpha \right) Y\left( t \right) +B_2\left( t,x,Y\left( t \right) \right) \right] dt+F_2\left( t, Y\left( t \right) \right) d\bar{W}^{Q_2}\left( t \right) \cr
&\quad+{\int_{\mathbb{Z}}{G_2}\left( t, Y\left( t \right) ,z \right)}{\tilde{N}_{{2}^{'}}}\left( dt,dz \right), \qquad\qquad\qquad\qquad Y\left( s \right)=y,
\end{align}
where 
$$
\bar{W}^{Q_2}\left( t \right) =\left\{ \begin{array}{l}
W_{1}^{Q_2}\left( t \right),  \\
W_{2}^{Q_2}\left( -t \right),  \\
\end{array} \right. \begin{array}{c}
\text{if}\ t\ge 0,\\
\text{if}\ t<0,\\
\end{array}
$$
$$
{\tilde{N}_{{2}^{'}}}\left( t,z \right) =\left\{ \begin{array}{l}
{\tilde{N}_{{1}^{'}}}\left( t,z \right), \\
{\tilde{N}_{{3}^{'}}}\left( -t,z \right),  \\
\end{array} \right. \begin{array}{c}
\text{if}\ t\ge 0,\\
\text{if}\ t<0,\\
\end{array}
$$
where ${\tilde{N}_{{1}^{'}}}\left( dt,dz \right)$ and ${\tilde{N}_{{3}^{'}}}\left( dt,dz \right)$ has the same  L\'{e}vy measure. The process $W_{1}^{Q_2}\left( t \right)$, $W_{2}^{Q_2}\left( t \right)$, ${\tilde{N}_{{1}^{'}}}\left( dt,dz \right)$ and  ${\tilde{N}_{{3}^{'}}}\left( dt,dz \right)$ are independent and the definition of which are given in Section \ref{sec-2}. According to the prove of \thmref{th3.2} in Section \ref{sec-3}, we can get that there exists a unique mild solution $ Y^x\left( \cdot;s,y\right) $ for equation (\ref{en41}) in the following form：
\begin{align}
Y^x\left( t;s,y \right) &= U_{\alpha }\left( t,s \right) y+\psi _{\alpha }\left( Y^x\left( \cdot;s,y \right) ;s \right)\left( t \right) +\int_s^t{U_{\alpha }\left( t,r \right) B_2\left( r,x,Y^x\left( r;s,y \right) \right)}dr\cr
&\quad+\int_s^t{U_{\alpha }\left( t,r \right) F_2\left( r,x,Y^x\left( r;s,y \right) \right)}d\bar{W}^{Q_2}\left( r \right) \cr
&\quad+\int_s^t{\int_{\mathbb{Z}}{U_{\alpha}\left( t,r \right) G_2\left( r,x,Y^x\left( r;s,y \right) ,z \right)}}\tilde{N}_{{2}^{'}}\left( dr,dz \right).\nonumber
\end{align}
Using the same argument as Lemma 4.1 in \cite{Xu2018Averaging}, we also can get that there exists $ \delta >0 $,  such that for any $ p\ge 1 $, have 
\begin{eqnarray}\label{en413}
\mathbb{E}\left\| Y^x\left( t;s,y \right) \right\|  ^{p}\le C_p\big( 1+\left\| x \right\| ^{p}+e^{-\delta p\left( t-s \right)}\left\| y \right\| ^{p} \big), \quad s<t.
\end{eqnarray}

Next, we giving the following auxiliary problem to prove that there exists an evolution family
of measures for the frozen fast equation: 
\begin{align}\label{en42}
dY\left( t \right) &=\left[ \left( A_2\left( t \right) -\alpha \right) Y\left( t \right) +B_2\left( t,x,Y\left( t \right) \right) \right] dt+F_2\left( t, Y\left( t \right) \right) d\bar{W}^{Q_2}\left( t \right) \cr
&\quad+\int_{\mathbb{Z}}{G_2}\left( t, Y\left( t \right) ,z \right) \tilde{N}_{{2}^{'}}\left( dt,dz \right),
\end{align}
for every $ s<t $ and $ t\in\mathbb{R}, $ we have
\begin{align}
Y^x\left( t \right) &=U_{\alpha   }\left( t,s \right) Y^x\left( s \right)+\psi _{\alpha  }\left( Y^x;s \right)\left( t \right) + \int_s^t{U_{\alpha  }\left( t,r \right) B_2\left( r,x ,Y^x\left( r \right) \right)}dr \cr
&\quad + \int_s^t{U_{\alpha   }\left( t,r \right) F_2\left( r, Y^x\left( r \right) \right)}d\bar{W}^{Q_2}\left( r \right) \cr
&\quad +\int_s^t{\int_{\mathbb{Z}}{U_{\alpha   }\left( t,r \right) G_2\left( r, Y^x\left( r \right) ,z \right)}}\tilde{N}_{{2}^{'}} \left( dr,dz \right).
\end{align}

Using the same arguments as  \cite{cerrai2017averaging,Xu2018Averaging}, we can establish the following crucial results for the frozen fast equation (\ref{en41}) whose proof will be presented in the Appendix. 
\begin{lem}\label{lem4.1}
	Under the assumptions {\rm (A1)-(A5)}, for any $ p\ge 1 $, there exists some $ \eta ^x\left( t \right) \in L^p\left( \varOmega ;E \right) $ such that
	\begin{eqnarray}\label{en414}
	\underset{s\rightarrow -\infty}{\lim}\mathbb{E}\left\| Y^x\left( t;s,y \right) -\eta ^x\left( t \right) \right\|^{p}=0,
	\end{eqnarray}
	where $ \eta ^x $ is a mild solution in $ \mathbb{R} $ of equation (\ref{en42}). 
	Moreover, there also exists some  $ \delta _p>0 $, such that
	\begin{eqnarray}\label{en415}
	\mathbb{E}\left\| Y^x\left( t;s,y \right) -\eta ^x\left( t \right) \right\|^{p}\le C_pe^{-\delta _p\left( t-s \right)}\left( 1+\left\| x \right\|^{p}+\left\| y \right\|^{p} \right), 
	\end{eqnarray} 
	and for any $ R>0,  $ there exists $ C_{R}>0 $ such that
	\begin{align}\label{419}
	\left\| x_{1}\right\| , \left\| x_{2}\right\|  \leq R\Longrightarrow\sup _{t \in \mathbb{R}} \mathbb{E}\left\| \eta^{x_{1}}(t)-\eta^{x_{2}}(t)\right\| ^{2} \leq C_{R}\left\| x_{1}-x_{2}\right\| ^{2}.
	\end{align}
\end{lem}

Now, for any  $ t\in \mathbb{R} $ and  $ x\in E$, we denote the law of the random variable  $ \eta ^x\left( t \right)  $ is $ \mu _{t}^x,  $  and introduce the transition evolution operator as
$$
P_{s,t}^{x}\varphi \left( y \right) =\mathbb{E}\varphi \left( Y^x\left( t;s,y \right) \right), \quad s<t, \ y\in E,
$$
where $ \varphi \in \mathcal{B}_b\left( E \right) $.
 
For any  $ p\ge 1 $, thanks to (\ref{en414}) and (\ref{en413}),  let $ s\rightarrow -\infty, $ we can get
\begin{eqnarray}\label{en426}
\underset{t\in \mathbb{R}}{\sup}\mathbb{E}\left\| \eta ^x\left( t \right) \right\| ^{p}\le C_p\left( 1+\left\| x \right\|^{p} \right),\quad x\in E,
\end{eqnarray}
therefore 
\begin{eqnarray}\label{en427}
\underset{t\in \mathbb{R}}{\sup}\int_{E}{\left\| y \right\|^{p}}\mu _{t}^{x}\left( dy \right) \le C_p\left( 1+\left\| x \right\|^{p} \right),\quad x\in E.
\end{eqnarray}

Then, under the assumptions (A1)-(A6), thanks to the \lemref{lem4.1} and the equation (\ref{en427}), 
the argument used in \cite{cerrai2017averaging} and \cite{Xu2018Averaging} can be adapted to the present situation, it is possible to prove that $\left\{ \mu _{t}^x  \right\}_{t \in \mathbb{R} }$ introduced above defines  an evolution family of probability measures for equation (\ref{en41}) and for any Lipschitz function $ \varphi  $ we also can get  
\begin{eqnarray}\label{en429}
\Big| P_{s,t}^{x}\varphi \left( y \right) -\int_{E}{\varphi \left(w \right)}\mu _{t}^{x}\left( dw \right) \Big|\le  L_\varphi e^{-\delta_1 \left( t-s \right)}\left( 1+\lVert x \rVert +\lVert y \rVert \right).
\end{eqnarray}
where $ L_{\varphi}=\underset{x\ne y}{ \sup }\frac{\left| \varphi \left( x \right) -\varphi \left( y \right) \right|}{\left| x-y \right|}. $ 
Moreover, 
the mapping  $t \mapsto \mu _{t}^{x}  \left(   t\in \mathbb{R},  x\in E \right)  $ 
is almost periodic and for any compact set  $ K\subset E, $ 
the family of functions
\begin{eqnarray}\label{427}
\left\{ t\in \mathbb{R}\mapsto \int_{E}{B_1\left( x,y \right)}\mu _{t}^{x}\left( dy \right): \ x\in K \right\} 
\end{eqnarray}
is uniformly almost periodic.


Therefore, according to the Theorem 3.4 in \cite{cerrai2017averaging}, we can define 
\begin{eqnarray}\label{en014} 
\bar{B}_1\left( x \right) :=\underset{T\rightarrow \infty}{\lim}\frac{1}{T}\int_0^T{\int_E{B_1\left(  x,y \right)}\mu _{t}^{x}\left( dy \right) dt, \quad x\in E},
\end{eqnarray}
and introduce the  averaged equation as follows:
\begin{eqnarray}\label{en67}
dX\left( t \right) &=&\left[ A_1(t)X\left( t \right) +\bar{B}_1\left( X\left( t \right) \right) \right] dt+F_1\left(t, X\left( t \right) \right) dW^{Q_1}\left( t \right) \cr
&&+\int_{\mathbb{Z}}{G_1\left( t,X\left( t \right) ,z \right)}\tilde{N}_1\left( dt,dz \right), \quad\quad  X\left( 0 \right) =x\in  E.
\end{eqnarray}
Due to the assumption {\rm (A3)} and  the equation (\ref{en427}), it is easy to prove 
\begin{eqnarray}\label{en62}
\left\| \bar{B}_1\left( x \right) \right\| \le c\left( 1+\left\| x \right\|^{m_1}\right).
\end{eqnarray} 
Moreover, we can also prove that the mapping $ \bar{B}:E \rightarrow E $ is locally Lipschitz continuous.
Indeed,  for any $ x_1, x_2 \in E, $ due to the assumption (A3) and (\ref{en426}), we have
\begin{align}
\bar{B}\left(x_{1}\right)-\bar{B}\left(x_{2}\right)&=\lim _{T \rightarrow \infty} \frac{1}{T} \int_{0}^{T} \mathbb{E}\left(B_{1}\left(x_{1}, \eta^{x_{1}}(t)\right)-B_{1}\left(x_{2}, \eta^{x_{2}}(t)\right)\right) dt\cr
&\quad \leq \lim _{T \rightarrow \infty} \frac{C}{T} \int_{0}^{T}  \mathbb{E}\Big[\left(\left\|x_{1}-x_{2}\right\|+\left\|\eta^{x_{1}}(t)-\eta^{x_{2}}(t)\right\|\right)\cr
&\qquad\qquad\qquad\qquad\times\big(1+\left\|\left( x_{1},\eta^{x_{1}}(t)\right) \right\|_{E\times E}^{\kappa}+\left\|\left( x_{2},\eta^{x_{2}}(t)\right) \right\|_{E\times E}^{\kappa}\big)\Big] dt\cr
&\leq\lim _{T \rightarrow \infty} \frac{C}{T} \int_{0}^{T}  \big(1+\left\|x_{1}\right\|^{\kappa}+\left\|x_{2}\right\|^{\kappa}\big)\big(\left\|x_{1}-x_{2}\right\|+\big(\sup _{t \in \mathbb{R}}\mathbb{E}\left\|\eta^{x_{1}}(t)-\eta^{x_{2}}(t)\right\|^{2}\big)^{\frac{1}{2}}\big)dt.\nonumber
\end{align}
Hence, thanks to (\ref{419}),  for any $ R > 0, $ it follows that 
\begin{align}\label{en432}
\left\| x_{1}\right\|, \left\| x_{2}\right\|  \leq R \Longrightarrow\left\| B_{1}\left(x_{1}\right)-B_{1}\left(x_{2}\right)\right\|  \leq C_{R}\left\| x_{1}-x_{2}\right\| .
\end{align}
Then, by proceeding as the Theorem \ref{th3.2}, we also can  prove that equation  (\ref{en67}) admits a unique mild solution $\bar{X} $ and it satisfy
\begin{eqnarray}\label{431}
\mathbb{E}\underset{t\in \left[ 0,T \right]}{\sup}\left\| \bar{X}\left( t \right) \right\|^{p}\le C_{p,T}，
\left( 1+\left\| x \right\|^{p}+\left\| y \right\|^{p} \right).
\end{eqnarray}


\subsection{Estimates of the auxiliary process $ \hat{Y}^\epsilon\left( t\right)  $}
Inspired by Khasminskii's idea in \cite{khas1968on}, for any $ \epsilon>0, $ we divide the interval $ [0,T] $ into subintervals of size $ \delta_\epsilon>0, $ where $ \delta_\epsilon $ is a fixed number depending on $ \epsilon. $ For each time interval  $ t\in \left[k\delta_\epsilon,   (k+1)\delta_\epsilon\land T \right], k=0,1,\cdots ,\lfloor  {T}/{\delta _{\epsilon}} \rfloor,  $ 
we construct the following auxiliary fast motion $ \hat{Y}^{\epsilon} $ 
with
initial value $ \hat{Y}^{\epsilon }\left( k\delta _{\epsilon} \right): = Y^{\epsilon  }\left( k\delta _{\epsilon} \right): $
\begin{eqnarray}\label{en72}
d\hat{Y}^{\epsilon }\left( t \right) &=&\frac{1}{\epsilon}\big[ \left( A_2\left( t \right) -\alpha \right) \hat{Y}^{\epsilon }\left( t \right) +B_{2 }\big( t, X^{\epsilon }\left( k\delta _{\epsilon} \right) ,\hat{Y}^{\epsilon }\left( t \right) \big) \big] dt\cr
&&+\frac{1}{\sqrt{\epsilon}}F_2\big( t,\hat{Y}^{\epsilon }\left( t \right) \big) dW ^{Q_2}\left( t \right)+\int_{\mathbb{Z}}{G_2}\big( t,\hat{Y}^{\epsilon }\left( t \right), z\big) \tilde{N}_{2}^{\epsilon}\left( dt,dz \right).
\end{eqnarray}
For any  $ p\ge 1, $ using the same argument as \thmref{th3.2},  we also can prove that  
\begin{eqnarray}\label{en73}
\int_0^T{\mathbb{E}\big\| \hat{Y}^{\epsilon }\left( t \right) \big\|^{p}}dt\le C_{p,T}\left( 1+\left\| x \right\|^{p}+\left\| y \right\|^{p} \right) .
\end{eqnarray}

Now, we prove the approximation result of $ Y^\epsilon\left( t\right) -\hat{Y}^\epsilon\left( t\right).   $
\begin{lem}\label{lem5.1}Under the assumptions (A1)-(A6), for any $ n\in \mathbb{N}, \ p\geq 1  $ and  any initial value $ x\in D\big(\left(-A_1 \right)^\theta  \big)     $  with $\theta\in [0,\bar{\theta}\land 1/p),   $ there exists some positive constant $C_{p, n,\theta, T}$, such that for any $\epsilon \in \left( 0,1 \right]$, we have 
	\begin{eqnarray}\label{en74}
	\mathbb{E}\big(\big\| \hat{Y}^{\epsilon }\left( t \right) - Y^{\epsilon }\left( t \right) \big\| ^{p}I_{\left\lbrace 0\leq t\leq T\land\tau_{n}^{\epsilon}\right\rbrace } \big) \le   C_{p, n,\theta, T} {\epsilon}^{-1}\delta _{\epsilon}^{1+\frac{p}{2}\land \theta p\land 1}
	 e^{C_{p  } \delta _{\epsilon}/\epsilon}\big( 1+\left\| x \right\|^{  m_1 p  }+\left\| y \right\|^{m_1 p } \big), 
	\end{eqnarray}
	where
\begin{align}
\tau _{n}^{\epsilon}:=\text{inf}\left\{ t\geq 0:\left\|  X^{\epsilon}\left( t \right) \right\| \geq n \right\} .\nonumber
\end{align}	
\end{lem}
\para{Proof:} To prove this, we need to construct some auxiliary processes.
For any $ n\in \mathbb{N} $ and $ \sigma_1,\sigma_2\in\mathbb{R}, $ we define 
\begin{eqnarray}
b_{1,n}\left(  \xi ,\sigma _1,\sigma _2 \right) :=\left\{ \begin{array}{c}
b_1\left(  \xi ,\sigma _1,\sigma _2 \right) ,\\
b_1\left(  \xi ,\sigma _1n/\left| \sigma _1 \right|,\sigma _2 \right) ,\\
\end{array} \right. \quad \begin{array}{c}
\text{if\ } \left| \sigma _1 \right|\leq n,\\
\text{if\ } \left| \sigma _1 \right|>n,\\
\end{array}
\end{eqnarray}
and
\begin{eqnarray}
b_{2,n}\left( t,\xi ,\sigma _1,\sigma _2 \right) :=\left\{ \begin{array}{c}
b_2\left( t,\xi ,\sigma _1,\sigma _2 \right) ,\\
b_2\left( t,\xi ,\sigma _1n/\left| \sigma _1 \right|,\sigma _2 \right) ,\\
\end{array} \right. \quad \begin{array}{c}
\text{if\ } \left| \sigma _1 \right|\leq n,\\
\text{if\ } \left| \sigma _1 \right|>n.\\
\end{array}
\end{eqnarray}
For each $ b_{i,n} $, denote the corresponding composition operator is $ B_{i,n}. $ Then, for any $ t\in \mathbb{R} $ and $  y\in E, $ we have
\begin{eqnarray}\label{en061}
\left\| x\right\| \leq n \Rightarrow 
B_{1,n}\left(  x,y \right) =B_1\left(  x,y \right) ,  \quad 
B_{2,n}\left( t,x,y \right) =B_2\left( t,x,y \right).
\end{eqnarray}
It is easy to get that the mapping $ B_{1,n} $ and $ B_{2,n} $ satisfy all conditions in (A3) and (A4), respectively. Moreover, for any fixed $  t   \in  \mathbb{R}   $ and $ \sigma_2\in \mathbb{R} $, the mapping $ B_{i,n}\left( t  ,\cdot,\sigma _2 \right) $ are Lipschitz-continuous.

Now, for any $ n\in\mathbb{N} $, we introduce the following system
\begin{eqnarray}\label{en0611}
\begin{split}
\begin{cases}
dX \left( t \right) &=\left[ A_1  X \left( t \right) +B_{1,n}\left(  X \left( t \right) ,Y \left( t \right) \right) \right] dt+F_{1}\left(  X \left( t \right) \right) dW ^{Q_1}\left( t \right)  +\int_{\mathbb{Z}}{G_{1}\left(  X \left( t \right) ,z \right)}\tilde{N}_1\left( dt,dz \right), \\
dY \left( t \right) &=\frac{1}{\epsilon}\left[ \left( A_2\left( t \right) -\alpha \right) Y \left( t \right) +B_{2,n}\left( t,X \left( t \right) ,Y \left( t \right) \right) \right] dt+\frac{1}{\sqrt{\epsilon}}F_2\left( t,Y \left( t \right) \right) dW ^{Q_2}\left( t \right) \\
&\quad   +\int_{\mathbb{Z}}{G_2}\left( t,Y \left( t \right)  ,z\right) \tilde{N}_{2}^{\epsilon}\left( dt,dz \right), \\
X \left( s \right) &=x, \quad Y \left( s \right) =y,\nonumber
\end{cases}
\end{split}
\end{eqnarray}
and denote the solution  is $ \left( X^{\epsilon,n}, Y^{\epsilon,n}\right)  $. 

Moreover, in each time interval  $ \left[ k\delta _{\epsilon},\left( k+1 \right) \delta _{\epsilon} \right] ,k=0,1,\cdots ,\lfloor  {T}/{\delta _{\epsilon}} \rfloor,  $ we define $ \hat{Y}^{\epsilon,n}   $ is the solution of the following problem  
\begin{align}\label{en75}
dY \left( t \right) &= \frac{1}{\epsilon}\left[ \left( A_2\left( t \right) -\alpha \right) Y\left( t \right) +B_{2,n }\left( t, X^{\epsilon,n }\left( k\delta _{\epsilon} \right) ,Y\left( t \right) \right) \right] dt+\frac{1}{\sqrt{\epsilon}}F_2\left( t,Y\left( t \right) \right) dW ^{Q_2}\left( t \right)\cr
&\quad+\int_{\mathbb{Z}}{G_2}\left( t,Y\left( t \right), z\right) \tilde{N}_{2}^{\epsilon}\left( dt,dz \right) , \qquad\qquad\qquad\qquad\qquad
Y\left( k\delta _{\epsilon} \right)  = Y^{\epsilon }\left( k\delta _{\epsilon} \right).
\end{align} 
Then, due to (\ref{en061}), it is easy to know that
\begin{align}\label{1}
&\quad\ \mathbb{E} \big\| \hat{Y}^{\epsilon,n }\left( t \right) - Y^{\epsilon,n }\left( t \right) \big\| ^{p}\cr &= \mathbb{E}\big(\big\| \hat{Y}^{\epsilon,n }\left( t \right) - Y^{\epsilon,n }\left( t \right) \big\| ^{p}I_{\left\lbrace 0\leq t\leq T\land\tau_{n}^{\epsilon}\right\rbrace } \big)+\mathbb{E}\big(\big\| \hat{Y}^{\epsilon,n }\left( t \right) - Y^{\epsilon,n }\left( t \right) \big\| ^{p}I_{\left\lbrace T\land\tau_{n}^{\epsilon}\leq t\leq T \right\rbrace } \big)\cr
&= \mathbb{E}\big(\big\| \hat{Y}^{\epsilon  }\left( t \right) - Y^{\epsilon  }\left( t \right) \big\| ^{p}I_{\left\lbrace 0\leq t\leq T\land\tau_{n}^{\epsilon}\right\rbrace } \big)+\mathbb{E}\big(\big\| \hat{Y}^{\epsilon,n }\left( t \right) - Y^{\epsilon,n }\left( t \right) \big\| ^{p}I_{\left\lbrace T\land\tau_{n}^{\epsilon}\leq t\leq T \right\rbrace } \big)\cr
&\geq \mathbb{E}\big(\big\| \hat{Y}^{\epsilon  }\left( t \right) - Y^{\epsilon  }\left( t \right) \big\| ^{p}I_{\left\lbrace 0\leq t\leq T\land\tau_{n}^{\epsilon}\right\rbrace } \big).  
\end{align}
So, we can prove the approximation result of $ Y^{\epsilon }\left( t\right) -\hat{Y}^{\epsilon }\left( t\right)    $ by geting an approximation result of $ Y^{\epsilon,n}\left( t\right) -\hat{Y}^{\epsilon,n}\left( t\right).    $

Fixed $ \epsilon >0, $  for any $ t\in \left[ k\delta _{\epsilon},  \left( k+1 \right) \delta _{\epsilon} \right], k=0,1,\cdots ,\lfloor {  T  }/{\delta _{\epsilon}} \rfloor $, let $ \rho_{\epsilon,n }\left( t \right)  $  be the solution of the following problem
\begin{align}
d\rho_{\epsilon,n }\left( t \right) &= \frac{1}{\epsilon}\left( A_2\left( t \right) -\alpha \right) \rho_{\epsilon,n }\left( t \right) dt+\frac{1}{\sqrt{\epsilon}}K_{\epsilon,n }\left( t \right) dW ^{Q_2}\left( t \right) +\int_{\mathbb{Z}}{H_{\epsilon,n }\left( t,z \right)}\tilde{N}_{2}^{\epsilon}\left( dt,dz \right), \quad  \rho_{\epsilon,n }\left( k\delta _{\epsilon} \right) =0,\nonumber
\end{align} 
where
$$
K_{\epsilon,n }\left( t \right) :=F_2\big( t, \hat{Y}^{\epsilon,n }\left( t \right) \big) -F_2\left( t, Y^{\epsilon,n }\left( t \right) \right), 
$$
$$
H_{\epsilon,n }\left( t,z \right) :=G_2\big( t, \hat{Y}^{\epsilon,n }\left( t \right) ,z \big) -G_2\left( t, Y^{\epsilon,n }\left( t \right) ,z \right). 
$$
We have 
$$
\rho _{\epsilon,n }\left( t \right) =\psi_{\alpha ,\epsilon }\left( \rho_{\epsilon,n };0 \right) \left( t \right) +\varGamma _{\epsilon,n }\left( t \right) +\varPsi _{\epsilon,n }\left( t \right), \quad t\in \left[ k\delta _{\epsilon},  \left( k+1\right) \delta _{\epsilon} \right], 
$$
where
$$
\varGamma _{\epsilon,n }\left( t \right) =\frac{1}{\sqrt{\epsilon}}\int_{0}^t{U_{\alpha ,\epsilon }\left( t,r \right) K_{\epsilon,n }\left( r \right)}dW^{Q_2}\left( r \right), 
$$
$$
\varPsi _{\epsilon,n }\left( t \right) =\int_{0}^t{\int_{\mathbb{Z}}{U_{\alpha ,\epsilon }\left( t,r \right) H_{\epsilon,n }\left( r,z \right)}}\tilde{N}_{2}^{\epsilon}\left( dr,dz \right). 
$$
For any $ s<r<t, $ we have $ U_{\alpha,\epsilon}\left( t,s\right) = U_{\alpha,\epsilon}\left( t,r\right) U_{\alpha,\epsilon}\left( r,s\right). $ Hence, the factorization arguement used in the the proof of Lemma 6.3  in \cite{cerrai2011averaging} and Lemma 6.2 in \cite{Xu2018Averaging} can also be used in present situation, and it is possible to show that 
\begin{eqnarray}\label{en202}
\mathbb{E}
\underset{r\in \left[ k\delta _{\epsilon},t \right]}{\text{sup}}
\left\| \varGamma _{\epsilon,n }\left(r \right) \right\| ^{p}+\mathbb{E}
\underset{r\in \left[ k\delta _{\epsilon},t \right]}{\text{sup}}\left\| \varPsi _{\epsilon,n }\left( t \right) \right\|^{p} \leq \frac{C_{p }}{\epsilon}\int_{k\delta _{\epsilon}}^{t}{\mathbb{E}\big\| \hat{Y}^{\epsilon,n }\left( r \right) - Y^{\epsilon,n }\left( r \right) \big\| ^{p} }dr.
\end{eqnarray}
Moreover, because $ \alpha $ is large enough, as a consequence of (\ref{56}) and  thanks to (\ref{en202}), it
holds that 
\begin{align}\label{57}
\mathbb{E}\underset{r\in \left[ k\delta _{\epsilon},t \right]}{\text{sup}}\left\| \rho _{\epsilon,n }\left( r\right) \right\|^p&\le \frac{C_{p }}{\epsilon}\int_{k\delta_{\epsilon}}^t{ \mathbb{E}\lVert Y^{\epsilon,n }\left( r \right) -\hat{Y}^{\epsilon,n }\left( r \right) \rVert ^p}dr.
\end{align}

If we denote $ \varLambda _{\epsilon,n }\left( t \right) :=\hat{Y}^{\epsilon,n }\left( t \right) - Y^{\epsilon,n }\left( t \right)  $  and $ \vartheta _{\epsilon,n }\left( t \right) := \varLambda _{\epsilon,n }\left( t \right) -  \rho _{\epsilon,n }\left( t \right)$, we have
\begin{align}
d\vartheta _{\epsilon,n}\left( t \right) &= \frac{1}{\epsilon}\big[ \left( A_2\left( t \right) -\alpha \right) \vartheta _{\epsilon,n }\left( t \right) +B_{2,n }\big( t, X^{\epsilon,n }\left( k\delta _{\epsilon} \right) ,\hat{Y}^{\epsilon,n }\left( t \right) \big)   \cr
&\qquad \quad   -B_{2,n }\left( t, X^{\epsilon,n }\left( t \right) , Y^{\epsilon,n }\left( t \right) \right) \big] dt, \qquad\qquad \vartheta _{\epsilon,n }\left( {k\delta_{\epsilon}} \right)=0.\nonumber
\end{align}
Then, it yields
\begin{align}
\frac{d}{dt}^-\lVert \vartheta _{\epsilon,n }\left( t \right) \rVert &=\frac{1}{\epsilon}\left< \left( A_2\left( t \right) -\alpha \right) \vartheta _{\epsilon,n }\left( t \right) ,\delta _{\vartheta _{\epsilon,n }\left( t \right)} \right>\cr
&\quad +\frac{1}{\epsilon}\big< B_{2,n }\big( t, X^{\epsilon,n }\left( k\delta _{\epsilon} \right) ,\hat{Y}^{\epsilon,n }\left( t \right) \big) -B_{2,n }\big( t, X^{\epsilon,n }\left( t \right) ,\hat{Y}^{\epsilon,n }\left( t \right) \big) ,\delta _{\vartheta _{\epsilon,n }\left( t \right)} \big> \cr
&\quad +\frac{1}{\epsilon}\big< B_{2,n }\big( t, X^{\epsilon,n }\left( t \right) ,\hat{Y}^{\epsilon,n }\left( t \right) \big)- B_{2,n }\big( t, X^{\epsilon,n }\left( t \right) , Y^{\epsilon,n }\left( t \right) \big)  ,\delta _{\vartheta _{\epsilon,n }\left( t \right)} \big>.\nonumber
\end{align}
For any $  t\in \left[ k\delta _{\epsilon}, \left( k+1 \right) \delta _{\epsilon} \land T \right], $ using the same argument as the proof of Lemma 6.2 in \cite{cerrai2011averaging},   due to (\ref{216}) and (\ref{en223}), we also can get
\begin{align}\label{en13}
\lVert \vartheta _{\epsilon,n }\left( t \right) \rVert&\le  \frac{C_n }{\epsilon}\int_{k\delta_{\epsilon}}^{t}{ e^{-\frac{\alpha}{\epsilon}\left( t-r \right)}\lVert X^{\epsilon,n }\left( k\delta _{\epsilon} \right)  - X^{\epsilon,n }\left( r \right)    \rVert }dr \cr
&\quad+\frac{1}{\epsilon}\int_{k\delta_{\epsilon}}^t{\exp \Big( -\frac{1}{\epsilon}\int_r^t{\tau _{\epsilon,n }\left( s \right)}ds \Big)}\tau _{\epsilon,n }\left( r \right) \lVert \rho _{\epsilon,n }\left( r \right) \rVert dr,
\end{align}
where
\begin{align}
\tau _{\epsilon,n }\left( t \right) :=\tau \big( t,\xi _{\epsilon,n }\left( t \right) , X^{\epsilon,n }\left( t,\xi _{\epsilon,n }\left( t \right) \right) ,\hat{Y}^{\epsilon,n }\left( t,\xi _{\epsilon,n }\left( t \right) \right) , Y^{\epsilon,n }\left( t,\xi _{\epsilon,n }\left( t \right) \right) \big),\nonumber
\end{align}
and $ \xi _{\epsilon,n }\left( t \right) $ is a point in $ \bar{\mathcal{O}} $ such that
\begin{align}
\left| \vartheta _{\epsilon,n }\left( t,\xi _{\epsilon,n }\left( t \right) \right) \right|=\lVert \vartheta _{\epsilon,n }\left( t \right) \rVert.\nonumber
\end{align}
By the Gronwall lemma and \lemref{lem3.4}, this implies
\begin{align}\label{en510}
\mathbb{E}\lVert \vartheta _{\epsilon,n }\left( t \right) \rVert ^p&\le \frac{C_{p,n}}{\epsilon ^p}\Big( \int_{k\delta _{\epsilon}}^t{e^{-\frac{\alpha}{\epsilon}\frac{p}{p-1}\left( t-r \right)}}dr \Big) ^{p-1}\int_{k\delta _{\epsilon}}^t{\mathbb{E}\lVert X^{\epsilon,n }\left( k\delta _{\epsilon} \right) - X^{\epsilon,n }\left( r \right) \rVert ^p}dr\cr
&\quad +\mathbb{E}\underset{r\in \left[ k\delta _{\epsilon},t \right]}{\text{sup}}\lVert \rho _{\epsilon,n }\left( r \right) \rVert ^p\Big( \frac{1}{\epsilon}\int_{k\delta_{\epsilon}}^t{\exp \Big( -\frac{1}{\epsilon}\int_r^t{\tau _{\epsilon,n }\left( s \right)}ds \Big)}\tau _{\epsilon,n }\left( r \right) dr \Big) ^p\cr
&\le \frac{C_{p,n,\theta,T}}{\epsilon}\delta _{\epsilon}^{1+\frac{p}{2}\land \theta p\land 1}\left( 1+\left\| x \right\|^{  m_1 p  }+\left\| y \right\|^{m_1 p } \right) +\frac{C_{p }}{\epsilon}\int_{k\delta _{\epsilon}}^t{\mathbb{E}\lVert \hat{Y}^{\epsilon,n }\left( r \right) - Y^{\epsilon,n }\left( r \right) \rVert ^p}dr.
\end{align}
According to the definition of $ \hat{Y}^{\epsilon,n}\left( t \right) - Y^{\epsilon,n}\left( t \right), $ thanks to (\ref{57}) and (\ref{en510}), we can get 
\begin{align}
\mathbb{E}\big\| \hat{Y}^{\epsilon,n}\left( t \right) - Y^{\epsilon,n}\left( t \right) \big\| ^{p} 
\leq  \frac{C_{p,n,\theta,T}}{\epsilon}\delta _{\epsilon}^{1+\frac{p}{2}\land \theta p\land 1}\left( 1+\left\| x \right\|^{  m_1 p  }+\left\| y \right\|^{m_1 p } \right) +\dfrac{C_{p }}{\epsilon }  \int_{k\delta _{\epsilon}}^{t}{\mathbb{E}\big\| \hat{Y}^{\epsilon,n}\left( r \right) - Y^{\epsilon,n}\left( r \right) \big\|^{p} }dr,\nonumber
\end{align}
using the Gronwall inequality and due to (\ref{1}), this implies the equation (\ref{en74}). \qed

\subsection{The proof of the main result}
In this part, the proof of our main result can be finished. That is, we will prove that the slow motion $ X^{\epsilon}  $ strongly converges to the averaged motion $ \bar{X}, $ as $ \epsilon \rightarrow 0 $.
\para{Proof of \thmref{th5.3}:}
For any $ p\geq 1  $ and $ n\in \mathbb{N}, $ we know that
\begin{align}
\mathbb{E}\Big( \underset{t\in \left[ 0,T \right]}{\text{sup}}\lVert X^{\epsilon}\left( t \right) -\bar{X}\left( t \right) \rVert ^p \Big) &\leq \mathbb{E}\Big( \underset{t\in \left[ 0,T\land \tilde{\tau} _{n}^{\epsilon} \right]}{\text{sup}}\lVert X^{\epsilon}\left( t \right) -\bar{X}\left( t \right) \rVert ^p \Big) +\mathbb{E}\Big( \underset{t\in \left[ 0,T \right]}{\text{sup}}\lVert X^{\epsilon}\left( t \right) -\bar{X}\left( t \right) \rVert ^p I_{\left\{ T>\tilde{\tau} _{n}^{\epsilon} \right\}} \Big),\nonumber 
\end{align}
where
\begin{align*}
\tilde{\tau} _{n}^{\epsilon}:=\text{inf}\left\{ t\geq 0:\left\|  X^{\epsilon}\left( t \right) \right\| +\left\|  \bar{X}\left( t \right) \right\|\geq n \right\} .
\end{align*}
Then, thanks to (\ref{en31}) and (\ref{431}), we can obtain
\begin{align} 
&\qquad\mathbb{E}\Big(\underset{t\in \left[ 0,T \right]}{\text{sup}}
\left\|  X^{\epsilon}\left( t \right) -\bar{X}\left( t \right) \right\| ^p I_{\left\lbrace T>\tilde{\tau} _{n}^{\epsilon}\right\rbrace } \Big)\cr 
&\le \Big[ \mathbb{E}\Big( \underset{t\in \left[ 0,T \right]}{\text{sup}}\left\|  X^{\epsilon}\left( t \right) -\bar{X}\left( t \right) \right\| ^{2p} \Big) \Big] ^{1/2}\left[ \mathbb{P}\left( T>\tilde{\tau} _{n}^{\epsilon} \right) \right] ^{1/2}\cr
&\le \frac{C_p}{n}\Big[ \mathbb{E}\Big(\underset{t\in \left[ 0,T \right]}{\text{sup}}\left\|  X^{\epsilon}\left( t \right) \right\| ^{2p}\Big)+\mathbb{E}\Big(\underset{t\in \left[ 0,T \right]}{\text{sup}}\left\|  \bar{X}\left( t \right) \right\| ^{2p}\Big) \Big] ^{1/2}\Big[ \mathbb{E}\Big(\underset{t\in \left[ 0,T \right]}{\text{sup}}\left\|  X^{\epsilon}\left( t \right) \right\|^2+\underset{t\in \left[ 0,T \right]}{\text{sup}}\left\|  \bar{X}\left( t \right) \right\|^2\Big) \Big] ^{1/2}\cr
&\le \frac{C_{p,T}}{n}\big( 1+\left\|  x \right\| ^{p+1}+\left\|  y \right\| ^{p+1} \big), 
\end{align}
let  $ n $ large enough, it yields
\begin{align}\label{en513}
\mathbb{E}\Big(\underset{t\in \left[ 0,T \right]}{\text{sup}}
\left\|  X^{\epsilon}\left( t \right) -\bar{X}\left( t \right) \right\| ^p I_{\left\lbrace T>\tilde{\tau} _{n}^{\epsilon}\right\rbrace } \Big)\rightarrow 0.
\end{align}

In order to prove \thmref{th5.3}, we shall prove the following result:
\begin{lem}\label{lem5.2}Under the assumptions {\rm (A1)-(A6)}, for any $ p\geq 1, $ we have 
	\begin{align}\label{en71}
	\underset{\epsilon \rightarrow 0}{\lim}\mathbb{E}\Big( \underset{t\in \left[ 0,T\land \tilde{\tau} _{n}^{\epsilon} \right]}{{\sup}}\left\|  X^{\epsilon}\left( t \right) -\bar{X}\left( t \right) \right\| ^p\Big) =0.
	\end{align}
\end{lem}
\para{Proof: } 
From the definition of $ X^\epsilon\left( t\right)  $ and $ \bar{X}\left( t\right),  $ we can get
\begin{align}\label{en515}
&\ \quad\mathbb{E}\Big( \underset{t\in \left[ 0,T\land \tilde{\tau} _{n}^{\epsilon} \right]}{{\sup}}\left\|  X^{\epsilon}\left( t \right) -\bar{X}\left( t \right) \right\| ^p\big)\cr
&\leq C_p\mathbb{E}\Big( \underset{t\in \left[ 0,T \land \tilde{\tau} _{n}^{\epsilon} \right]}{{\sup}} \Big\|  \int_0^{t}{e^{A_1 \left( t -r \right)}\left( B_1\left(   X^{\epsilon}\left( r \right) , Y^{\epsilon}\left( r \right) \right) -\bar{B}_1\left(   \bar{X}\left( r \right) \right) \right)}dr \Big\| ^p\Big)\cr
&\quad+C_p\mathbb{E}\Big( \underset{t\in \left[ 0,T\land \tilde{\tau} _{n}^{\epsilon} \right]}{{\sup}} \Big\|  \int_0^{t}{e^{A_1 \left( t -r \right)}\left( F_1\left(   X^{\epsilon}\left( r \right)  \right) -F_1\left(   \bar{X}\left( r \right) \right) \right)}dW^{Q_1}\left( r \right) \Big\| ^p\Big)\cr
&\quad+C_p\mathbb{E}\Big( \underset{t\in \left[ 0,T\land \tilde{\tau} _{n}^{\epsilon}  \right]}{{\sup}} \Big\|  \int_0^{t}{\int_\mathbb{Z}{}e^{A_1 \left( t -r \right)}\left( G_1\left(   X^{\epsilon}\left( r \right) ,z \right) -G_1\left(   \bar{X}\left( r \right),z \right) \right)}\tilde{N}_1\left( dr,dz \right) \Big\| ^p\Big)\cr
&:= \mathcal{K}_1\left( t\right) +\mathcal{K}_2\left( t\right) +\mathcal{K}_3\left( t\right).
\end{align} 
Step 1: We estimate the term $ \mathcal{K}_{1}\left( t \right), $ we know that
\begin{align}\label{en511}
\mathcal{K}_1\left( t\right)&\le C_p\mathbb{E}\Big( \underset{t\in \left[ 0,T\land \tilde{\tau} _{n}^{\epsilon} \right]}{{\sup}}   \Big\|  \int_0^t{e^{A_1 \left( t-r \right)}\big( B_1\left(   X^{\epsilon}\left( r \right) , Y^{\epsilon}\left( r \right) \right) -B_1\big(   X^{\epsilon}\left( r\left( \delta _{\epsilon} \right) \right) ,\hat{Y}^{\epsilon}\left( r \right) \big) \big)}dr \Big\|^p\Big) \cr
& \quad+C_p\mathbb{E}\Big( \underset{t\in \left[ 0,T\land \tilde{\tau} _{n}^{\epsilon} \right]}{{\sup}} \Big\|  \int_0^t{e^{A_1 \left( t-r \right)}\big( B_1\big(   X^{\epsilon}\left( r\left( \delta _{\epsilon} \right) \right) ,\hat{Y}^{\epsilon}\left( r \right) \big) -\bar{B}_1\left( X^{\epsilon}\left( r\left( \delta _{\epsilon} \right) \right) \right) \big)}dr \Big\| ^p\Big)\cr
& \quad+C_p \mathbb{E}\Big( \underset{t\in \left[ 0,T\land \tilde{\tau} _{n}^{\epsilon} \right]}{{\sup}}\Big\|  \int_0^t{e^{A_1 \left( t-r \right)}\left( \bar{B}_1\left(  X^{\epsilon}\left( r\left( \delta _{\epsilon} \right) \right) \right) -\bar{B}_1\left(  \bar{X}\left( r \right) \right) \right)}dr \Big\| ^p\Big)\cr
&:=C_p\sum_{i=1}^3{\mathbb{E}\Big( \underset{t\in \left[ 0,T\land \tilde{\tau} _{n}^{\epsilon} \right]}{{\sup}}\left\| \mathcal{I}_i\left( t \right)\right\|^p \Big)}. 
\end{align}
where $ r(\delta_{\epsilon})=\lfloor r/\delta _{\epsilon} \rfloor \delta _{\epsilon} $ and $ \lfloor r  \rfloor $ denotes the largest integer which is no more than $ r. $  \\
For $ \mathcal{I}_1\left( t \right), $ using the H\"{o}lder inequality and in view of (\ref{en0221}), it is possible to get that
\begin{align}\label{511} 
\mathbb{E}\Big( \underset{t\in \left[ 0,T\land \tilde{\tau} _{n}^{\epsilon} \right]}{{\sup}} \left\| \mathcal{I}_1\left( t \right)\right\|^p\Big)  
&\leq  C_{p, T}\mathbb{E}\Big(\int_0^{T\land \tilde{\tau} _{n}^{\epsilon}}{\big\|   B_1\left(   X^{\epsilon}\left( r \right) , Y^{\epsilon}\left( r \right) \right) -B_1\big(  X^{\epsilon}\left( r\left( \delta _{\epsilon} \right) \right) ,\hat{Y}^{\epsilon}\left( r \right) \big)   \big\|^p} dr\Big)  \cr
&\le C_{p, T}\mathbb{E}\Big(\int_0^{T \land\tilde{\tau} _{n}^{\epsilon}  }{ \big( \left\|  X^{\epsilon}\left( r \right) - X^{\epsilon}\left( r\left( \delta _{\epsilon} \right) \right) \right\| ^p+ \big\|  Y^{\epsilon}\left( r \right) -\hat{Y}^{\epsilon}\left( r \right) \big\| ^p \big)}\cr
&\quad\qquad\qquad  \times \big( 1+ \left\|  \left( X^{\epsilon}\left( r \right), Y^{\epsilon}\left( r \right)\right)  \right\|_{E\times E} ^{p\kappa }+ \big\| \big( X^{\epsilon}\left( r\left( \delta _{\epsilon} \right) \right), \hat{Y}^{\epsilon}\left( r \right) \big) \big\|_{E\times E}^{p\kappa  }  \big) dr \Big)\cr 
&\le  C_{p, T}\int_0^{T }{\mathbb{E}\Big[\big( \left\|  X^{\epsilon}\left( r \right) - X^{\epsilon}\left( r\left( \delta _{\epsilon} \right) \right) \right\| ^p+ \big\|  Y^{\epsilon}\left( r \right) -\hat{Y}^{\epsilon}\left( r \right) \big\| ^p I_{\left\lbrace 0\leq r\leq T\land\tilde{\tau}_{n}^{\epsilon}\right\rbrace } \big)}\cr
&\quad\qquad\qquad  \times \big( 1+ \left\|  \left( X^{\epsilon}\left( r \right), Y^{\epsilon}\left( r \right)\right)  \right\|_{E\times E} ^{p\kappa }+ \big\| \big( X^{\epsilon}\left( r\left( \delta _{\epsilon} \right) \right), \hat{Y}^{\epsilon}\left( r \right) \big) \big\|_{E\times E}^{p\kappa  }  \big)\Big] dr \cr 
&\le  C_{p,n,\theta,T}   \big(\delta _{\epsilon}^{\frac{p}{2}\land p\theta \land 1}+{\epsilon}^{-1}\delta _{\epsilon}^{1+\frac{p}{2}\land \theta p\land 1}
e^{C_{p  } \delta _{\epsilon}/\epsilon}\big)\big( 1+\lVert x \rVert ^{pm_1+p\kappa}+\lVert y \rVert ^{pm_1+p\kappa} \big).
\end{align}
%
For $ \mathcal{I}_2\left( t \right), $ we can get
\begin{align}\label{en518}
\mathbb{E}\Big( \sup_{t\in \left[ 0,T\land \tilde{\tau} _{n}^{\epsilon} \right]}\lVert \mathcal{I}_2\left( t \right) \rVert ^p \Big) \le \Big[ \mathbb{E}\Big( \sup_{t\in \left[ 0,T\land \tilde{\tau} _{n}^{\epsilon} \right]}\lVert \mathcal{I}_2\left( t \right) \rVert ^{2p-2} \Big) \mathbb{E}\Big( \sup_{t\in \left[ 0,T\land \tilde{\tau} _{n}^{\epsilon} \right]}\lVert \mathcal{I}_2\left( t \right) \rVert ^2 \Big) \Big] ^{\frac{1}{2}}.
\end{align}
Thanks to (\ref{en31}), (\ref{431}) and (\ref{en73}), we obtain
\begin{eqnarray}\label{en519}
\mathbb{E}\Big( \sup_{t\in \left[ 0,T\land \tilde{\tau} _{n}^{\epsilon} \right]}\lVert \mathcal{I}_2\left( t \right) \rVert ^{2p-2} \Big) &\le& C_T\mathbb{E}\Big( \int_0^{T\land \tilde{\tau} _{n}^{\epsilon}}{\lVert B_1\big( X^{\epsilon}\left( r\left( \delta _{\epsilon} \right) \right) ,\hat{Y}^{\epsilon}\left( r \right) \big) -\bar{B}_1\left( X^{\epsilon}\left( r\left( \delta _{\epsilon} \right) \right) \right) \rVert ^{2p-2}}dr \Big)\cr
&\le& C_T  \int_0^T{\mathbb{E}\lVert B_1\big( X^{\epsilon}\left( r\left( \delta _{\epsilon} \right) \right) ,\hat{Y}^{\epsilon}\left( r \right) \big) \rVert ^{2p-2}+\mathbb{E}\lVert \bar{B}_1\left( X^{\epsilon}\left( r\left( \delta _{\epsilon} \right) \right) \right) \rVert ^{2p-2}}dr \cr
&\le& C_T\int_0^T{}\Big( 1+\mathbb{E}\underset{\sigma \in \left[ 0,r \right]}{\text{sup}}\lVert X^{\epsilon}\left( \sigma \right) \rVert ^{2pm_1-2m_1}+\mathbb{E}\lVert \hat{Y}^{\epsilon}\left( r \right) \rVert ^{2p-2}\Big)dr\cr 
&\le& C_{p,T}\left( 1+\lVert x \rVert ^{2pm_1 }+\lVert y \rVert ^{2pm_1 } \right). 
\end{eqnarray}
Moreover, we also have
\begin{align} 
\left\| \mathcal{I}_2\left( t \right) \right\| ^2
&\le   2   \Big\| \sum_{k=0}^{\lfloor t/\delta _{\epsilon} \rfloor-1}{ }e^{A_1\left( t-\left( k+1 \right) \delta _{\epsilon} \right)}
 \int_{k\delta _{\epsilon}}^{\left( k+1 \right) \delta _{\epsilon}}{e^{A_1\left( \left( k+1 \right) \delta _{\epsilon}-r \right)}
 	\big( B_1\big(  X^{\epsilon}\left( k\delta _{\epsilon} \right) ,\hat{Y}^{\epsilon}\left( r \right) \big) -\bar{B}_1\left( X^{\epsilon}\left( k\delta _{\epsilon} \right) \right) \big)}dr \Big\| ^2\cr
&\quad  +2 \Big\|  \int_{ t\left( \delta _{\epsilon}\right) }^{t}{e^{A_1 \left( t-r \right)}\big( B_1\big(  X^{\epsilon}\left( r\left( \delta _{\epsilon} \right) \right) ,\hat{Y}^{\epsilon}\left( r \right) \big) -\bar{B}_1\left( X^{\epsilon}\left( r\left( \delta _{\epsilon}\right)  \right) \right) \big)}dr \Big\| ^2\cr
&:=\mathcal{I}_{21}\left( t \right)+\mathcal{I}_{22}\left( t \right).\nonumber 
\end{align}
For the term $ \mathcal{I}_{21}\left( t \right), $ we can get
\begin{align}\label{512} 
&\qquad\mathbb{E}\Big( \underset{t\in \left[ 0,T\land \tilde{\tau} _{n}^{\epsilon} \right]}{{\sup}}\mathcal{I}_{21}\left( t \right)\Big)
\cr
&\leq C\lfloor\frac{T }{\delta _{\epsilon}}\rfloor\mathbb{E}\Big(  \sum_{k=0}^{\lfloor   \left( T\land \tilde{\tau} _{n}^{\epsilon}\right)    /\delta _{\epsilon} \rfloor-1}{ } \Big\|  \int_{k\delta _{\epsilon}}^{\left( k+1 \right) \delta _{\epsilon}}{e^{A_1 \left( \left( k+1 \right) \delta _{\epsilon}-r \right)}\big( B_1\big(  X^{\epsilon}\left( k\delta _{\epsilon} \right) ,\hat{Y}^{\epsilon}\left( r \right) \big) -\bar{B}_1\left( X^{\epsilon}\left( k\delta _{\epsilon} \right) \right) \big)}dr \Big\| ^2 \Big) \cr 
&\le   C \frac{T }{\delta _{\epsilon}} \mathbb{E}\Big( \sum_{k=0}^{\lfloor   \left( T\land \tilde{\tau} _{n}^{\epsilon}\right)   /\delta _{\epsilon} \rfloor-1}{}\Big\|  \int_0^{\delta _{\epsilon}}{e^{A_1 \left(  \delta _{\epsilon}-r \right)}\big( B_1\big(   X^{\epsilon}\left( k\delta _{\epsilon} \right) ,\hat{Y}^{\epsilon}\left( k\delta _{\epsilon}+r \right) \big)  -\bar{B}_1\left( X^{\epsilon}\left( k\delta _{\epsilon} \right) \right) \big) } dr \Big\|   ^2\Big)  \cr 
&\le   C_T\frac{\epsilon ^2 }{\delta _{\epsilon}}  \mathbb{E}\Big(\sum_{k=0}^{\lfloor  \left( T\land \tilde{\tau} _{n}^{\epsilon}\right) /\delta _{\epsilon} \rfloor-1}{} \Big\|  \int_0^{{\delta_{\epsilon}}/\epsilon}{e^{A_1 \left(  \delta _{\epsilon}-\epsilon r \right)}\big( B_1\big(  X^{\epsilon}\left( k\delta _{\epsilon} \right) ,\hat{Y}^{\epsilon}\left( k\delta _{\epsilon}+\epsilon r \right) \big)  -\bar{B}_1\left( X^{\epsilon}\left( k\delta _{\epsilon} \right) \right) \big)} dr \Big\|  ^2 \Big) \cr 
&\le C_T   \frac{\epsilon ^2}{\delta _{\epsilon}^2}  \mathbb{E}\Big(\underset{0\le k\le \lfloor   \left( T\land \tilde{\tau} _{n}^{\epsilon}\right)   /\delta _{\epsilon} \rfloor-1}{\max} \Big\| \int_0^{{\delta_{\epsilon}}/\epsilon}{e^{A_1 \left(  \delta _{\epsilon}-\epsilon r \right)}\Big( B_1\big(   X^{\epsilon}\left( k\delta _{\epsilon} \right) ,\hat{Y}^{\epsilon}\left( k\delta _{\epsilon}+\epsilon r \right) \big)  } \cr
& \quad\qquad \quad 
-\int_E{B_1\left( X^{\epsilon}\left( k\delta _{\epsilon} \right) ,y \right)}\mu _{r}^{ X^{\epsilon}\left( k\delta _{\epsilon} \right)}\left( dy \right)+ \int_E{B_1\left( X^{\epsilon}\left( k\delta _{\epsilon} \right) ,y \right)}\mu _{r}^{ X^{\epsilon}\left( k\delta _{\epsilon} \right)}\left( dy \right) 
-\bar{B}_1\left( X^{\epsilon}\left( k\delta _{\epsilon} \right) \right) \Big) dr \Big\|  ^2 \Big) \cr 
&\le  C_T\frac{\epsilon ^2}{\delta _{\epsilon}^{2}} \mathbb{E}\bigg[\underset{0\le k\le \lfloor  \left( T\land \tilde{\tau} _{n}^{\epsilon}\right)  /\delta _{\epsilon} \rfloor-1}{\max} \Big(   \Big\| \int_0^{\delta _{\epsilon}/\epsilon}{} e^{A_1\left(  \delta _{\epsilon}-\epsilon r \right)}   \Big( B_1\big( X^{\epsilon}\left( k\delta _{\epsilon} \right) ,\tilde{Y}^{ X^{\epsilon}\left( k\delta _{\epsilon} \right) }\left( r;0 , Y^{\epsilon}\left( k\delta _{\epsilon} \right) \right) \big)   \cr
& \qquad\qquad\qquad\qquad\qquad\qquad\qquad\qquad\qquad \qquad\qquad\qquad -\int_E{B_1\left( X^{\epsilon}\left( k\delta _{\epsilon} \right) ,y \right)}\mu _{r}^{ X^{\epsilon}\left( k\delta _{\epsilon} \right)}\left( dy \right) 
\Big) dr \Big\|  ^2  \cr
& \qquad\qquad\qquad\qquad+   \Big\| \int_0^{\delta _{\epsilon}/\epsilon}{e^{A_1\left(  \delta _{\epsilon}-\epsilon r \right)}\Big( \int_E{B_1\left( X^{\epsilon}\left( k\delta _{\epsilon} \right) ,y \right)}\mu _{r}^{ X^{\epsilon}\left( k\delta _{\epsilon} \right)}\left( dy \right) -\bar{B}_1\left( X^{\epsilon}\left( k\delta _{\epsilon} \right) \right) \Big)}dr \Big\| ^2\Big)  \bigg]\cr
&\le  C_T\mathbb{E}\Big(\underset{0\le k\le \lfloor\left(  T\land \tilde{\tau} _{n}^{\epsilon}\right) /\delta _{\epsilon} \rfloor-1}{\max} \Big\|\frac{1  }{\delta_{\epsilon}/\epsilon } \int_0^{\delta _{\epsilon}/\epsilon}{ \Big( \int_E{B_1\left( X^{\epsilon}\left( k\delta _{\epsilon} \right) ,y \right)}\mu _{r}^{ X^{\epsilon}\left( k\delta _{\epsilon} \right)}\left( dy \right) -\bar{B}_1\left( X^{\epsilon}\left( k\delta _{\epsilon} \right) \right) \Big)}dr \Big\| ^2 \Big) \cr
&\quad+  C_T\frac{\epsilon ^2}{\delta _{\epsilon}^{2}}
\underset{0\le k\le \lfloor   T  /\delta _{\epsilon} \rfloor-1}{\max}\int_0^{\delta _{\epsilon}/\epsilon}{\int_r^{\delta _{\epsilon}/\epsilon}{\mathcal{J}_k\left(  \sigma,r \right)}d\sigma dr}, 
\end{align}
where 
\begin{align}
\mathcal{J}_k\left(  \sigma,r  \right) &=\mathbb{E}  \Big[ e^{A_1\left(  \delta _{\epsilon}-\epsilon r \right)}\Big( B_1\big( X^{\epsilon}\left( k\delta _{\epsilon} \right),\tilde{Y}^{ X^{\epsilon}\left( k\delta _{\epsilon} \right)}\left( r;s, Y^{\epsilon}\left( k\delta _{\epsilon} \right) \right) \big) -\int_E{B_1\left( X^{\epsilon}\left( k\delta _{\epsilon} \right),w \right)}\mu _{r}^{ X^{\epsilon}\left( k\delta _{\epsilon} \right)}\left( dw \right) \Big)   \cr
&\quad\times e^{A_1\left(  \delta _{\epsilon}-\epsilon \sigma \right)}\Big( B_1\big( X^{\epsilon}\left( k\delta _{\epsilon} \right),\tilde{Y}^{ X^{\epsilon}\left( k\delta _{\epsilon} \right)}\left( \sigma;s, Y^{\epsilon}\left( k\delta _{\epsilon} \right) \right) \big) -\int_E{B_1\left( X^{\epsilon}\left( k\delta _{\epsilon} \right),w \right)}\mu _{\sigma}^{ X^{\epsilon}\left( k\delta _{\epsilon} \right)}\left( dw \right) \Big) \Big].\nonumber
\end{align}
and $ \tilde{Y}^{ X^{\epsilon}\left( k\delta _{\epsilon} \right) }\left( r;0 , Y^{\epsilon}\left( k\delta _{\epsilon} \right) \right) $ is the solution of the fast motion equation (\ref{en41}) with the initial datum given by  $ Y^{\epsilon}\left( k\delta _{\epsilon} \right) $ and the frozen slow component given by $ X^{\epsilon}\left( k\delta _{\epsilon} \right), $ the distribution of it coincides with the distribution of $ \hat{Y}^{\epsilon}\left( k\delta _{\epsilon}+\epsilon r \right). $ 

Then, we can adapt the proof of appendix A in \cite{Fu2011Strong} to the present situation, and it is possible to show that
\begin{align}
\mathcal{J}_k\left(  \sigma,r \right)&\le c e^{-\frac{\delta _2}{2}\left( \sigma-r\right) } \big( 1+\mathbb{E}\left\|  X^{\epsilon}\left( k\delta _{\epsilon} \right) \right\|^{2\kappa\lor 2\lor m_1} +\mathbb{E}\left\|  Y^{\epsilon}\left( k\delta _{\epsilon} \right) \right\|^{2\kappa\lor 2} \big)^2\cr
&\le C_{T} e^{-\frac{\delta _2}{2}\left( \sigma-r\right) } \big( 1+\left\|x \right\|^{4\kappa\lor 4\lor 2 m_1} +\left\|y \right\|^{4\kappa\lor 4\lor 2 m_1} \big),\nonumber
\end{align}
it follows that
\begin{align}\label{516}
\mathbb{E}\Big( \underset{t\in \left[ 0,T\land \tilde{\tau} _{n}^{\epsilon} \right]}{{\sup}}\mathcal{I}_{21}\left( t \right)\Big)\leq C_T\frac{\epsilon}{\delta_\epsilon}\big( 1+\left\|x \right\|^{4\kappa\lor 4\lor 2 m_1} +\left\|y \right\|^{4\kappa\lor 4\lor 2 m_1} \big).
\end{align}
Moreover, thanks to (\ref{en31}) and (\ref{en73}), we can get
\begin{eqnarray}\label{523}
\mathbb{E}\Big( \underset{t\in \left[ 0,T\land \tilde{\tau} _{n}^{\epsilon} \right]}{{\sup}}\mathcal{I}_{22}\left( t \right)\Big)
&\leq & C\delta_{\epsilon}\mathbb{E}\Big( \underset{t\in \left[ 0,T\land \tilde{\tau} _{n}^{\epsilon} \right]}{{\sup}}   \int_{  t\left( \delta _{\epsilon}\right) }^{t}{\big\|  B_1\big(  X^{\epsilon}\left( r\left( \delta _{\epsilon}\right)  \right) ,\hat{Y}^{\epsilon}\left( r \right) \big) \big\|^2+\big\| \bar{B}_1\left( X^{\epsilon}\left( r\left( \delta _{\epsilon}\right) \right) \right)  \big\| ^2}dr\Big)\cr
& \leq & C\delta_{\epsilon}\int_{0}^{T}{\Big(1+\mathbb{E}  \underset{\sigma\in \left[ 0,r \right]}{{\sup}}\big\| X^\epsilon\left( \sigma\right) \big\|^{2m_1}+\mathbb{E}\big\|\hat{Y}^\epsilon\left( r\right) \big\|^2\Big)}dr\cr
& \leq&  C_T\delta_{\epsilon} \big( 1+\left\| x \right\|^{2m_1}+\left\| y \right\|^{2m_1} \big).
\end{eqnarray}
For the term $ \mathcal{I}_{3}\left( t \right), $ according to the (\ref{en432}) and \lemref{lem3.4}, we have
\begin{align}\label{518}
&\qquad\mathbb{E}\Big( \underset{t\in \left[ 0,T\land \tilde{\tau} _{n}^{\epsilon} \right]}{{\sup}} \left\| \mathcal{I}_3\left( t \right)\right\|^p\Big)  \cr
&\le C_p \mathbb{E}\Big( \underset{t\in \left[ 0,T\land \tilde{\tau} _{n}^{\epsilon} \right]}{{\sup}}\Big\|  \int_0^t{e^{A_1 \left( t-r \right)}\left( \bar{B}_1\left( X^{\epsilon}\left( r\left( \delta _{\epsilon} \right) \right) \right) -\bar{B}_1\left(  X^{\epsilon}\left( r \right) \right) \right)}dr \Big\| ^p\Big)\cr
&\quad+C_p \mathbb{E}  \Big(\underset{t\in \left[ 0,T\land \tilde{\tau} _{n}^{\epsilon} \right]}{{\sup}}\Big\|  \int_0^t{e^{A_1 \left( t-r \right)}\left( \bar{B}_1\left(  X^{\epsilon}\left( r \right) \right) -\bar{B}_1\left(  \bar{X}\left( r \right) \right) \right)}dr \Big\| ^p\Big) \cr
&\le C_{p, T} \mathbb{E} \int_0^{T\land \tilde{\tau} _{n}^{\epsilon}}{\left\|  \bar{B}_1\left(  X^{\epsilon}\left( r\left( \delta _{\epsilon} \right) \right) \right) -\bar{B}_1\left(  X^{\epsilon}\left( r \right) \right) \right\| ^p}dr   +C_{ p,T} \mathbb{E} \int_0^{T\land \tilde{\tau} _{n}^{\epsilon}}{\left\|  \bar{B}_1\left(  X^{\epsilon}\left( r \right) \right) -\bar{B}_1\left(  \bar{X}\left( r \right) \right) \right\| ^p}dr \cr
&\le  C_{ p,T} \int_0^{T }{\mathbb{E}\Big(\underset{\sigma\in \left[ 0,r\land \tilde{\tau} _{n}^{\epsilon} \right]}{{\sup}}\left\|  X^{\epsilon}\left( \sigma \right) -\bar{X}\left( \sigma \right) \right\| ^p\Big)}dr +C_{ p,T}\delta _{\epsilon}^{ {p}/{2}\land p\theta\land 1 }\big( 1+\left\| x\right\|^{pm_1}+\left\| y\right\|^{pm_1}\big). 
\end{align}
Step 2: We estimate the terms $ \mathcal{K}_{2}\left( t \right)  $ and $ \mathcal{K}_{3}\left( t \right). $ For $ \mathcal{K}_{2}\left( t \right), $ we can get
\begin{align}
\mathcal{K}_{2}\left( t \right)&=C_p\mathbb{E}\Big( \underset{t\in \left[ 0,T\land \tilde{\tau} _{n}^{\epsilon} \right]}{{\sup}} \Big\|  \int_0^{t}{e^{A_1 \left( t -r \right)}\left( F_1\left(   X^{\epsilon}\left( r \right)  \right) -F_1\left(   \bar{X}\left( r \right) \right) \right)}dW^{Q_1}\left( r \right) \Big\| ^p\Big)\cr
&=C_p\mathbb{E}\Big( \underset{t\in \left[ 0,T \right]}{{\sup}} \Big\|  \int_0^{t\land \tilde{\tau} _{n}^{\epsilon}}{e^{A_1 \left( t\land \tilde{\tau} _{n}^{\epsilon} -r \right)}\left( F_1\left(   X^{\epsilon}\left( r \right)  \right) -F_1\left(   \bar{X}\left( r \right) \right) \right)}dW^{Q_1}\left( r \right) \Big\| ^p\Big)\cr
&:=C_p\mathbb{E}\Big( \underset{t\in \left[ 0,T \right]}{{\sup}} \big\| \varGamma _{\epsilon}\left( t \right)   \big\| ^p\Big).\nonumber
\end{align}
Then, using a factorization argument for $ \varGamma _{\epsilon}\left( t \right), $ for  any $ \theta\in \left( 0,1/2\right),  $ we have
\begin{align}
\varGamma _{\epsilon}\left( t \right) =C_{\theta}\int_0^{t\land \tilde{\tau}_{n}^{\epsilon}}{\left( t\land \tilde{\tau}_{n}^{\epsilon}-r \right) ^{\theta -1}e^{A_1\left( t\land \tilde{\tau}_{n}^{\epsilon}-r \right)}\varLambda _{\epsilon ,\theta}\left( r \right)}dr,\nonumber
\end{align}
where
\begin{align*}
\varLambda _{\epsilon ,\theta}\left( r \right) :=\int_0^r{\left( r-s \right) ^{-\theta}e^{A_1\left( r-s \right)}\left( F_1\left( X^{\epsilon}\left( s \right) \right) -F_1\left( \bar{X}\left( s \right) \right) \right)}dW^{Q_1}\left( s \right).
\end{align*}
For any $ p> 1/\theta, $ we have
\begin{align}
&\quad\ \underset{t\in \left[ 0,T \right]}{\text{sup}}\lVert \varGamma _{\epsilon}\left( t \right) \rVert ^p\le C_{ p,\theta}\underset{t\in \left[ 0,T \right]}{\text{sup}}\Big[ \Big( \int_0^{t\land \tilde{\tau}_{n}^{\epsilon}}{\left( t\land \tilde{\tau}_{n}^{\epsilon}-r \right) ^{\frac{p\left( \theta -1 \right)}{p-1}}}dr \Big) ^{p-1}\int_0^{t\land \tilde{\tau}_{n}^{\epsilon}}{\lVert \varLambda _{\epsilon ,\theta}\left( r \right) \rVert ^p}dr \Big] \cr
&\le C_{p,\theta,T}\underset{t\in \left[ 0,T \right]}{\text{sup}}\int_0^t{\lVert \varLambda _{\epsilon ,\theta}\left( r \right) \rVert ^p I_{\left\{ r\le \tilde{\tau}_{n}^{\epsilon} \right\}}}dr\cr
&=C_{p,\theta,T}\underset{t\in \left[ 0,T \right]}{\text{sup}}\int_0^t{\Big\| \int_0^{r\land \tilde{\tau}_{n}^{\epsilon}}{\left( r\land \tilde{\tau}_{n}^{\epsilon}-s \right) ^{-\theta}e^{A_1\left( r\land \tilde{\tau}_{n}^{\epsilon}-s \right)}\left( F_1\left( X^{\epsilon}\left( s \right) \right) -F_1\left( \bar{X}\left( s \right) \right) \right)}dW^{Q_1}\left( s \right) \Big\| ^p}dr\cr
&=C_{p,\theta,T}\int_0^{\text{T}}{\Big\| \int_0^r{\left( r\land \tilde{\tau}_{n}^{\epsilon}-s \right) ^{-\theta}e^{A_1\left( r\land \tilde{\tau}_{n}^{\epsilon}-s \right)}\left( F_1\left( X^{\epsilon}\left( s \right) \right) -F_1\left( \bar{X}\left( s \right) \right) \right) I_{\left\{ \text{s}\le \tilde{\tau}_{n}^{\epsilon} \right\}}}dW^{Q_1}\left( s \right) \Big\| ^p}dr.\nonumber
\end{align}
Then, if we choose $ \bar{\theta}, $ such that
$ 2\bar{\theta} +\frac{\beta _1\left( \rho _1-2 \right)}{\rho _1}<1, $
we know that there exist some $ \bar{p}>1 $, such that  
$$ \frac{\bar{2\theta \bar{p}}}{\bar{p}-2}+\frac{\beta _1\left( \rho _1-2 \right)}{\rho _1}\frac{\bar{p}}{\bar{p}-2}<1. $$ According to the Burkholder-Davis-Gundy inequality and the equation (3.5) in \cite{cerrai2009khasminskii}, for any  $p\ge \bar{p},$ we have
\begin{eqnarray}
\mathcal{K}_2\left( t \right) &\le& C_{\theta ,p,T}\int_0^{\text{T}}{\mathbb{E}\Big( \int_0^r{\left( r\land \tilde{\tau}_{n}^{\epsilon}-s \right) ^{-2\theta}\lVert e^{A_1\left( r\land \tilde{\tau}_{n}^{\epsilon}-s \right)}\left( F_1\left( X^{\epsilon}\left( s \right) \right) -F_1\left( \bar{X}\left( s \right) \right) \right) \rVert _{2}^{2} I_{\left\{ \text{s}\le \tilde{\tau}_{n}^{\epsilon} \right\}}}ds \Big) ^{\frac{p}{2}}}dr\cr
&\le& C_{\theta ,p,T}\int_0^{\text{T}}{\mathbb{E}\Big( \int_0^r{\left( r\land \tilde{\tau}_{n}^{\epsilon}-s \right) ^{-2\theta -\frac{\beta _1\left( \rho _1-2 \right)}{\rho _1}}\lVert F_1\left( X^{\epsilon}\left( s \right) \right) -F_1\left( \bar{X}\left( s \right) \right) \rVert ^2 I_{\left\{ \text{s}\le \tilde{\tau}_{n}^{\epsilon} \right\}}}ds \Big) ^{\frac{p}{2}}}dr\cr
&\le& C_{\theta ,p,T}\int_0^{\text{T}}{\mathbb{E}\Big[ \Big( \int_0^{r\land \tilde{\tau}_{n}^{\epsilon}}{\left( r\land \tilde{\tau}_{n}^{\epsilon}-s \right) ^{-\frac{2\theta p}{p-2}-\frac{\beta _1\left( \rho _1-2 \right)}{\rho _1}\frac{p}{p-2}}}ds \Big) ^{\frac{p-2}{2}}}\cr
&&\qquad\qquad\qquad\qquad\times\int_0^r{\underset{\sigma \in \left[ 0,s\land \tilde{\tau}_{n}^{\epsilon} \right]}{\text{sup}}\lVert F_1\left( X^{\epsilon}\left( \sigma \right) \right) -F_1\left( \bar{X}\left( \sigma \right) \right) \rVert ^p}ds \Big]dr\cr
&\le& C_{\theta ,p,T}\int_0^{\text{T}}{\mathbb{E}\Big(\underset{\sigma \in \left[ 0,s\land \tilde{\tau}_{n}^{\epsilon} \right]}{\text{sup}}\lVert  X^{\epsilon}\left( \sigma \right)   -  \bar{X}\left( \sigma \right)   \rVert ^p\Big) ds},
\end{eqnarray}
where the last equation using Fubini’s theorem.

Moreover, using Kunita's first inequality for $\mathcal{K}_{3}\left( t \right), $ it is easy to get that
\begin{align}\label{0515}
 \mathcal{K}_3\left( t\right)
&\le  C_p\mathbb{E} \Big( \int_0^{T\land \tilde{\tau} _{n}^{\epsilon}}{\int_{\mathbb{Z}}{\lVert \left( G_1\left( X^{\epsilon}\left( r \right) ,z \right) -G_1\left( \bar{X}\left( r \right) ,z \right) \right) \rVert ^2}} v_1\left( dz \right) dr \Big) ^{\frac{p}{2}}  \cr
&\quad+C_p\mathbb{E}  \int_0^{T\land \tilde{\tau} _{n}^{\epsilon}}{\int_{\mathbb{Z}}{\lVert \left( G_1\left( X^{\epsilon}\left( r \right) ,z \right) -G_1\left( \bar{X}\left( r \right) ,z \right) \right) \rVert ^p}} v_1\left( dz \right) dr   \cr
&\le C_{p,T}\int_0^T{\mathbb{E}\Big(\sup_{\sigma \in \left[ 0,r\land \tilde{\tau} _{n}^{\epsilon} \right]}\lVert X^{\epsilon}\left( \sigma \right) -\bar{X}\left( \sigma \right) \rVert ^p\Big)}dr
\end{align}
Step 3: Combining (\ref{en511})-(\ref{0515}) and using the Gronwall inequality for (\ref{en515}), it yields
\begin{align}\label{519}
&\qquad\mathbb{E}\Big(\underset{0\le t\le T\land \tilde{\tau} _{n}^{\epsilon}}{\text{sup}}\left\|  X^{\epsilon}\left( t \right) -\bar{X}\left( t \right) \right\| ^p\Big)\cr
&\leq C_{p,n,\theta,T}   \big(\delta _{\epsilon}^{p/2\land p\theta \land 1}+{\epsilon}^{-1}\delta _{\epsilon}^{1+\frac{p}{2}\land \theta p\land 1}
e^{C_{p  } \delta _{\epsilon}/\epsilon}\big)\big( 1+\lVert x \rVert ^{pm_1+p\kappa}+\lVert y \rVert ^{pm_1+p\kappa} \big)\cr
&\quad +C_{p,T}\left( 1+\lVert x \rVert ^{pm_1 }+\lVert y \rVert ^{pm_1 } \right)\bigg[\frac{\epsilon}{\delta_\epsilon}\big( 1+\left\|x \right\|^{4\kappa\lor 4\lor 2 m_1} +\left\|y \right\|^{4\kappa\lor 4\lor 2 m_1} \big)+  \delta_{\epsilon} \big( 1+\left\| x \right\|^{2m_1}+\left\| y \right\|^{2m_1} \big) \cr
&\quad +  \mathbb{E}\Big(\underset{0\le k\le \lfloor\left(  T\land \tilde{\tau} _{n}^{\epsilon}\right) /\delta _{\epsilon} \rfloor-1}{\max} \Big\|\frac{1  }{\delta_{\epsilon}/\epsilon } \int_0^{\delta _{\epsilon}/\epsilon}{ \Big( \int_E{B_1\left( X^{\epsilon}\left( k\delta _{\epsilon} \right) ,y \right)}\mu _{r}^{ X^{\epsilon}\left( k\delta _{\epsilon} \right)}\left( dy \right) -\bar{B}_1\left( X^{\epsilon}\left( k\delta _{\epsilon} \right) \right) \Big)}dr \Big\| ^2 \Big)\bigg]^{\frac{1}{2}}\cr
&\quad+ C_{ p,T}\delta _{\epsilon}^{ {p}/{2}\land p\theta\land 1 }\big( 1+\left\| x\right\|^{pm_1}+\left\| y\right\|^{pm_1}\big).
\end{align}
Selecting  $ \delta _{\epsilon}=\epsilon \ln ^{\epsilon ^{-\kappa}} $, then if we take  $ \kappa <\frac{  {p}/{2}\land \theta p\land 1 }{ C_{p  }+ 1+ {p}/{2}\land \theta p\land 1 }, $ we have
\begin{align}\label{520}
\underset{\epsilon \rightarrow 0}{\lim}\big({\epsilon}^{-1}\delta _{\epsilon}^{1+\frac{p}{2}\land \theta p\land 1}
e^{C_{p  } \delta _{\epsilon}/\epsilon}+\epsilon/\delta_{\epsilon}\big)=0.
\end{align}
Moreover, according to Theorem 3.4 in \cite{cerrai2017averaging}, thanks to the family of functions (\ref{427}) is uniformly almost periodic, we can get that the limit
$$
\frac{1  }{T } \int_s^{s+T}{   \int_E{B_1\left(x,y \right)}\mu _{r}^{x}\left( dy \right)  }dr 
$$
converges to  $ \bar{B}_1\left( x \right)  $ uniformly with respect to $ s\in \mathbb{R}  $ and  $ x $ in any compact set $ K \subset E.  $ From the definition of $ \tilde{\tau} _{n}^{\epsilon}, $ we know that $ \left\|  X^{\epsilon}\left( k\delta _{\epsilon} \right) \right\| \leq n $ for any $ 0\le k\le \lfloor\left(  T\land \tilde{\tau} _{n}^{\epsilon}\right) /\delta _{\epsilon} \rfloor-1. $ Hence, we can get
\begin{eqnarray}\label{521}
\underset{\epsilon \rightarrow 0}{\lim}
\mathbb{E}\big(\underset{0\le k\le \lfloor\left(  T\land \tilde{\tau} _{n}^{\epsilon}\right) /\delta _{\epsilon} \rfloor-1}{\max}\Big\|\frac{1  }{\delta_{\epsilon}/\epsilon } \int_0^{\delta _{\epsilon}/\epsilon}{   \int_E{B_1\left( X^{\epsilon}\left( k\delta _{\epsilon} \right) ,y \right)}\mu _{r}^{ X^{\epsilon}\left( k\delta _{\epsilon} \right)}\left( dy \right)   }dr-\bar{B}_1\left( X^{\epsilon}\left( k\delta _{\epsilon} \right) \right) \Big\| ^2\big)=0.
\end{eqnarray}
Thanks to (\ref{519})-(\ref{521}), it follows that
\begin{align} 
\underset{\epsilon \rightarrow 0}{\lim}\mathbb{E}\Big(\underset{t\in \left[ 0,T\land \tilde{\tau} _{n}^{\epsilon} \right]}{\text{sup}}\left\|  X^{\epsilon}\left( t \right) -\bar{X}\left( t \right) \right\| ^p\Big)=0.\nonumber
\end{align}
The proof is complete.\qed

Now, letting $ \epsilon \rightarrow 0 $ firstly and $ n  \rightarrow \infty $ secondly, according to the equation (\ref{en513}) and (\ref{en71}), it is easy to get that (\ref{510}) holds.  This completes the proof of \thmref{th5.3}.\qed

\section*{Acknowledgments}
This work was partly supported by the National Natural Science Foundation of China under Grant No. 11772255, the Fundamental Research Funds for the Central Universities, the Research Funds for Interdisciplinary Subject of Northwestern Polytechnical University, the Shaanxi Project for Distinguished Young Scholars, the Shaanxi Provincial Key R\&D Program 2020KW-013 and 2019TD-010.

\section*{Appendix}
In this section, we give the detailed proofs of \lemref{lem3.3} and \lemref{lem4.1}:
\para{Proof of \lemref{lem3.3}:}
\setcounter{equation}{0}
\renewcommand{\theequation}{A\arabic{equation}}
For any $ t \in \left[ 0,T\right]  $ and $ \epsilon \in \left( 0,1 \right] , $ we denote 
$$
\varGamma _{1,\epsilon}\left( t \right) :=\int_0^t{e^{A_1\left( t-r \right)}F_1\left(r, X^{\epsilon,n}\left( r \right) \right)}dW^{Q_1}\left( r \right), 
$$
and
$$
\varPsi _{1,\epsilon}\left( t \right) :=\int_0^t{\int_{\mathbb{Z}}{e^{A_1\left( t-r \right)}G_1\left( r, X^{\epsilon,n}\left( r \right) ,z \right)}}\tilde{N}_1\left( dr,dz \right). 
$$
Set $\varLambda _{1,\epsilon}\left( t \right) := X^{\epsilon,n}\left( t \right) -\varGamma _{1,\epsilon}\left( t \right) -\varPsi _{1,\epsilon}\left( t \right)$, we have
\begin{align}
\frac{d}{dt}\varLambda _{1,\epsilon}\left( t \right) &=  A_1\varLambda _{1,\epsilon}\left( t \right) +B_{1,n}\left( \varLambda _{1,\epsilon}\left( t \right) +\varGamma _{1,\epsilon}\left( t \right) +\varPsi _{1,\epsilon}\left( t \right) , Y^{\epsilon,n}\left( t \right) \right), \quad \varLambda _{1,\epsilon}\left( 0 \right) =x. \nonumber
\end{align}
Due to (\ref{en02}) and (\ref{en07}), we can get
\begin{align}\label{en333}
\frac{d}{dt}^-\left\| \varLambda _{1,\epsilon}\left( t \right) \right\| 
&=  \left<  A_1\varLambda _{1,\epsilon}\left( t \right) ,\delta_{\varLambda _{1,\epsilon}\left( t \right)} \right> + \left< B_{1,n}\left(  \varLambda _{1,\epsilon}\left( t \right) +\varGamma _{1,\epsilon}\left( t \right) +\varPsi _{1,\epsilon}\left( t \right) , Y^{\epsilon,n}\left( t \right) \right) \right. \cr
&\quad\qquad \left. -B_{1,n}\left(  \varGamma _{1,\epsilon}\left( t \right) +\varPsi _{1,\epsilon}\left( t \right) , Y^{\epsilon,n}\left( t \right) \right) ,\delta_{\varLambda _{1,\epsilon}\left( t \right)} \right> \cr
&\quad+\left< B_{1,n}\left(  \varGamma _{1,\epsilon}\left( t \right) +\varPsi _{1,\epsilon}\left( t \right) , Y^{\epsilon,n}\left( t \right) \right),\delta_{\varLambda _{1,\epsilon}\left( t \right)}\right\rangle  \cr
&\leq  C \left\| \varLambda _{1,\epsilon}\left( t \right) \right\|   +c \left( 1+\left\| \varGamma _{1,\epsilon}\left( t \right) \right\|    +\left\| \varPsi _{1,\epsilon}\left( t \right) \right\|  +\left\| Y^{\epsilon,n}\left( t \right) \right\|  \right)\cr
&\quad+C\left\|  B_{1,n}\left(  \varGamma _{1,\epsilon}\left( t \right) +\varPsi _{1,\epsilon}\left( t \right) , Y^{\epsilon,n}\left( t \right) \right)\right\| \cr
&\leq C \left\| \varLambda _{1,\epsilon}\left( t \right) \right\|   +c \left( 1+\left\| \varGamma _{1,\epsilon}\left( t \right) \right\|^{m_1}   +\left\| \varPsi _{1,\epsilon}\left( t \right) \right\|^{m_1} +\left\| Y^{\epsilon,n}\left( t \right) \right\|  \right). 
\end{align}
Due to the Gronwall inequality, it yields
$$
\left\| \varLambda _{1,\epsilon}\left( t \right) \right\|  \le e^{C t}\left\| x \right\|  +C \int_0^t{e^{C \left( t-r \right)}\left( 1+\left\| \varGamma _{1,\epsilon}\left( r \right) \right\|^{m_1}  +\left\| \varPsi _{1,\epsilon}\left( r \right) \right\|^{m_1}  +\left\| Y^{\epsilon,n}\left( r \right) \right\|   \right)}dr.
$$
For any $ p\geq 1, $ using the H\"{o}lder inequality, we can get 
$$
\left\| \varLambda _{1,\epsilon}\left( t \right) \right\|^p  \le C_{p,T}\Big( 1+ \left\| x \right\|^p+\underset{t\in \left[ 0,T \right]}{\text{sup}} \left\| \varGamma _{1,\epsilon}\left( t \right) \right\|^{m_1p} +\underset{t\in \left[ 0,T \right]}{\text{sup}} \left\| \varPsi_{1,\epsilon}\left( t \right) \right\|^{m_1p}+\int_0^T{  \left\| Y^{\epsilon,n}\left( r \right) \right\|^p    }dr \Big).
$$
This implies that 
\begin{align}\label{en37}
\mathbb{E}\underset{t\in \left[ 0,T \right]}{\sup}\left\| X^{\epsilon,n}\left( t \right) \right\|^p
&\leq  C_{p,T } \Big( 1+\left\| x \right\|^p  +\mathbb{E}\underset{t\in \left[ 0,T \right]}{\sup}\left\| \varGamma _{1,\epsilon}\left( t \right) \right\|^{m_1p}  +\mathbb{E}\underset{t\in \left[ 0,T \right]}{\sup}\left\| \varPsi _{1,\epsilon}\left( t \right) \right\|^{m_1p}\cr
&\quad\qquad\quad+  \int_0^T{ \mathbb{E}\left\| Y^{\epsilon,n}\left( r \right) \right\|^p    }dr  \Big).
\end{align}
Under the assumption (A5), by proceeding as the Lemma 4.1 in \cite{cerrai2009khasminskii} and Lemma 3.1 in \cite{Xu2018Averaging}. It is possible to prove that for any $ p\geq 1, $ we have 
\begin{align}\label{en36}
\mathbb{E}\underset{t\in \left[ 0,T \right]}{\sup}\left\| \varGamma _{1,\epsilon}\left( t \right) \right\| ^{p}+\mathbb{E}\underset{t\in \left[ 0,T \right]}{\sup}\left\| \varPsi _{1,\epsilon}\left( t \right) \right\|  ^{p} \le  C_{p ,T } \int_0^T{ \left( 1+\mathbb{E}\lVert X^{\epsilon ,n}\left( r \right) \rVert ^{\frac{p}{m_1}}\right)  }dr.  
\end{align}
Substituting (\ref{en36}) into (\ref{en37}), we can get 
\begin{eqnarray}\label{en6}
\mathbb{E}\underset{t\in \left[ 0,T \right]}{\sup}\left\| X^{\epsilon,n}\left( t \right) \right\| ^{p}
\leq C_{p,T  }  \Big( 1+\left\| x \right\|  ^{p}+ \int_0^T{ \mathbb{E}\left\| Y^{\epsilon,n}\left( r \right) \right\|^p    }dr   \Big) +C_{p,T }   \int_0^T{ \mathbb{E}\underset{\sigma\in \left[ 0,r \right]}{\sup}\lVert X^{\epsilon ,n}\left( \sigma \right) \rVert ^p }dr.
\end{eqnarray}

Now, we estimate $  \int_0^T{ \mathbb{E}\left\| Y^{\epsilon,n}\left( r \right) \right\|^p    }dr. $ For any $ t\in \left[ 0,T\right] , $ we define
$$
\varGamma _{2,\epsilon}\left( t \right) :=\frac{1}{\sqrt{\epsilon}}\int_0^t{U_{\alpha ,\epsilon }\left( t,r \right) F_2\left( r, Y^{\epsilon,n}\left( r \right) \right) dW^{Q_2}\left( r \right)},
$$
and
$$
\varPsi _{2,\epsilon}\left( t \right) :=\int_0^t{\int_{\mathbb{Z}}{U_{\alpha ,\epsilon }\left( t,r \right) G_2\left( r, Y^{\epsilon,n}\left( r \right) ,z \right)}}\tilde{N}_{2}^{\epsilon}\left( dr,dz \right).
$$
As before, we set $ \varLambda _{2,\epsilon}\left( t \right) := Y^{\epsilon,n}\left( t \right) -\varGamma _{2,\epsilon}\left( t \right) -\varPsi _{2,\epsilon}\left( t \right), $ we have
\begin{align}
\frac{d}{dt}\varLambda _{2,\epsilon}\left( t \right) &= \frac{1}{\epsilon}\left( \gamma  \left( t \right) A_2-\alpha \right) \varLambda _{2,\epsilon}\left( t \right) +\frac{1}{\epsilon} L  \left( t \right) \left( \varLambda _{2,\epsilon}\left( t \right) +\varGamma _{2,\epsilon}\left( t \right) +\varPsi _{2,\epsilon}\left( t \right) \right) \cr
&\quad+\frac{1}{\epsilon}B_{2,n}\left( t, X^{\epsilon,n}\left( t \right) ,\varLambda _{2,\epsilon}\left( t \right) +\varGamma _{2,\epsilon}\left( t \right) +\varPsi _{2,\epsilon}\left( t \right) \right), \qquad\ \varLambda _{2,\epsilon}\left( 0 \right) =y.\nonumber
\end{align}
For any $ p\geq 1, $ by proceeding as the proof of  (\ref{en333}), thanks to the constant $ \alpha $ is large enpugh, according to (\ref{en02}) and (\ref{en06}),   we can get 
\begin{eqnarray} 
\frac{d}{dt}^-\left\| \varLambda _{2,\epsilon}\left( t \right) \right\|^p 
\leq  -\frac{\alpha p}{2\epsilon} \left\| \varLambda _{1,\epsilon}\left( t \right) \right\|^p   +\frac{C_p}{\epsilon} \left( 1+\left\| X^{\epsilon,n}\left( t \right) \right\|^p+\left\| \varGamma _{2,\epsilon}\left( t \right) \right\|^{m_2 p}   +\left\| \varPsi _{2,\epsilon}\left( t \right) \right\|^{m_2 p}   \right), 
\end{eqnarray} 
by comparison, this implies
\begin{eqnarray}
\lVert \varLambda _{2,\epsilon}\left( t \right) \rVert ^p\leq e^{-\frac{\alpha p}{2\epsilon}t}\lVert y \rVert ^p+\frac{C_p}{\epsilon}\int_0^t{e^{-\frac{\alpha p}{2\epsilon}\left( t-r \right)}\left( 1+\lVert X^{\epsilon ,n}\left( r \right) \rVert ^p+\lVert \varGamma _{2,\epsilon}\left( r \right) \rVert ^{m_2p}+\lVert \varPsi _{2,\epsilon}\left( r \right) \rVert ^{m_2p} \right)}dr.
\end{eqnarray}
Then, integrating both sides in time and using Young's inequality, it follows that
\begin{align}\label{en4}
\int_0^t{\lVert \varLambda _{2,\epsilon}\left( r \right) \rVert ^p}dr&\le \frac{C_p}{\epsilon}\int_0^t{\left( 1+\lVert X^{\epsilon ,n}\left( r \right) \rVert ^p+\lVert \varGamma _{2,\epsilon}\left( r \right) \rVert ^{m_2p}+\lVert \varPsi _{2,\epsilon}\left( r \right) \rVert ^{m_2p} \right)}dr\int_0^t{e^{-\frac{\alpha p}{2\epsilon}r}}dr\cr
&\quad+\int_o^t{e^{-\frac{\alpha p}{2\epsilon}r}\lVert y \rVert ^p}dr\cr
&\le C_{p,1}\left( t \right) \int_0^t{\Big( \underset{\sigma\in \left[ 0,r \right]}{\sup}\lVert X^{\epsilon ,n}\left( \sigma \right) \rVert ^p+\lVert \varGamma _{2,\epsilon}\left( r \right) \rVert ^{m_2p}+\lVert \varPsi _{2,\epsilon}\left( r \right) \rVert ^{m_2p} \Big)}dr \cr
&\quad +C_{p,T}\left( 1+\lVert y \rVert ^p \right),
\end{align}
and we can easy to know that $ C_{p,1}\left( t\right) $ is a continuous increasing function and $ C_{p,1}\left( 0 \right)=0. $

For any $ p\geq 1, $ adapt the proof of Proposition 4.2 in \cite{cerrai2009khasminskii} and Lemma 3.1 in \cite{Xu2018Averaging} to the present situation,  it is possible to get that
\begin{eqnarray}\label{en3}
\int_0^t{\mathbb{E}\lVert \varGamma _{2,\epsilon}\left( r \right) \rVert ^p}dr +\int_0^t{\mathbb{E}\lVert \varPsi _{2,\epsilon}\left( r \right) \rVert ^p}dr\le C_{p,2}\left( t\right) \int_0^t{\big( 1+\mathbb{E}\lVert Y^{\epsilon ,n}\left( r \right) \rVert ^{\frac{p}{m_2}} \big)}dr,
\end{eqnarray}
where $ C_{p,2}\left( t\right) $ is a continuous increasing function and $ C_{p,2}\left( 0 \right)=0. $ 
Hence, thanks to the equation (\ref{en4}) and (\ref{en3}), it yields
\begin{align}
\int_0^t{\mathbb{E}\lVert  Y^{ \epsilon,n}\left( r \right) \rVert ^p}dr&\leq C_p\int_0^t{\mathbb{E}\lVert \varLambda _{2,\epsilon}\left( r \right) \rVert ^p}dr+C_p\int_0^t{\mathbb{E}\lVert \varGamma _{2,\epsilon}\left( r \right) \rVert ^p}dr+C_p\int_0^t{\mathbb{E}\lVert \varPsi _{2,\epsilon}\left( r \right) \rVert ^p}dr\cr
&\leq C_{p,T}\left( 1+\lVert y \rVert ^p \right) +C_{p,3}\left( t \right) \int_0^t{\big(\mathbb{E}\underset{\sigma\in \left[ 0,r \right]}{\sup}\lVert X^{\epsilon ,n}\left( \sigma \right) \rVert ^p+\mathbb{E}\lVert  Y^{ \epsilon,n}\left( r \right) \rVert ^{p}  \big)}dr.\nonumber
\end{align}
As $ C_{p,3}\left( t \right) $ is continuous increasing function and vanishes at $ t=0 ,$ we can fix $t_1>0$, such that for any $t\le t_1$, have $ C_{p,3}\left( t \right) \le {1}/{2} $. Then
\begin{align}\label{en5}
\int_0^t{\mathbb{E}\lVert  Y^{ \epsilon,n}\left( r \right) \rVert ^p}dr \leq C_{p,T}\left( 1+\lVert y \rVert ^p \right) +  \int_0^t{ \mathbb{E}\underset{\sigma\in \left[ 0,r \right]}{\sup}\lVert X^{\epsilon ,n}\left( \sigma \right) \rVert ^p }dr, \quad t\in \left[ 0,t_1 \right].
\end{align}
Substituting (\ref{en5}) into (\ref{en6}), we can get
\begin{eqnarray} 
\mathbb{E}\underset{r\in \left[ 0,t \right]}{\sup}\left\| X^{\epsilon,n}\left( r \right) \right\|  ^{p}\le C_{p,T} \left( 1+\left\| x \right\| ^{p}+\left\| y \right\| ^{p} \right)+C_{p,T }   \int_0^t{ \mathbb{E}\underset{\sigma\in \left[ 0,r \right]}{\sup}\lVert X^{\epsilon ,n}\left( \sigma \right) \rVert ^p }dr, \quad t\in \left[ 0,t_1 \right],\nonumber 
\end{eqnarray}
so
\begin{eqnarray}\label{en7}
\mathbb{E}\underset{r\in \left[ 0,t \right]}{\sup}\left\| X^{\epsilon,n}\left( r \right) \right\| ^{p}\le C_{p,T} \left( 1+\left\| x \right\| ^{p}+\left\| y \right\| ^{p} \right), \quad t\in \left[ 0,t_1 \right]. 
\end{eqnarray}
Using this for (\ref{en5}), we have
\begin{align} 
\int_0^t{\mathbb{E}\lVert  Y^{ \epsilon,n}\left( r \right) \rVert ^p}dr \le C_{p,T} \left( 1+\left\| x \right\| ^{p}+\left\| y \right\|  ^{p} \right), \quad t\in \left[ 0,t_1 \right]. 
\end{align}
Repeating this proof process in the intervals $\left[ t_1,2t_1 \right] ,  \left[ 2t_1,3t_1 \right] $ etc., we can get (\ref{en33}) and (\ref{en34}) hold. 
The proof of \lemref{lem3.3} is complete. \qed
\para{Proof of \lemref{lem4.1}:} Fix  $ h>0 $ and define
$$
\rho \left( t  \right) =Y^x\left( t;s,y \right) -Y^x\left( t;s-h,y \right), \quad s<t.
$$
It is easy to know that  $ \rho \left( t \right) $ is the unique mild solution of the following problem
\begin{eqnarray}\label{en416} 
\begin{split}
\begin{cases}
d\rho \left( t \right) &=\left[ \left( A_2\left( t \right) -\alpha \right) \rho \left( t \right) -J^x\left( t \right) \rho \left( t \right) \right] dt+K^x\left( t \right) \rho \left( t \right) d\bar{W}^{Q_2}\left( t \right) +\int_{\mathbb{Z}}{H^x\left( t,z \right) \rho \left( t \right)}\tilde{N}_{{2}^{'}}\left( dt,dz \right)\\
\rho \left( s \right) &=y-Y^x\left( s;s-h,y \right),
\end{cases} 
\end{split}\nonumber
\end{eqnarray}
where 
\begin{align}
J^x\left( t,\xi \right) &=\frac{b_2\left( t,\xi ,x\left( \xi \right) ,Y^x\left( t;s,y \right)\left( \xi \right) \right) -b_2\left( t,\xi ,x\left( \xi \right),Y^x\left( t;s-h,y \right)\left( \xi \right) \right)}{\rho \left( t \right) \left( \xi \right)}\cr
&= \tau  \left( t,\xi ,x\left( \xi \right) ,Y^x\left( t;s,y \right) \left( \xi \right) ,Y^x\left( t;s-h,y \right) \left( \xi \right) \right), \quad \xi \in \mathcal{O},
\end{align}
$$
K^x\left( t,\xi \right) =\frac{f_2\left( t,\xi , Y^x\left( t;s,y \right) \left( \xi \right) \right) -f_2\left( t,\xi , Y^x\left( t;s-h,y \right) \left( \xi \right) \right)}{\rho \left( t \right) \left( \xi \right)},\quad \xi \in \mathcal{O},
$$
$$
H^x\left( t,\xi,z \right) =\frac{g_2\left( t,\xi , Y^x\left( t;s,y \right) \left( \xi \right) ,z \right) -g_2\left( t,\xi , Y^x\left( t;s-h,y \right) \left( \xi \right) ,z \right)}{\rho \left( t \right) \left( \xi \right)}, \quad z \in \mathbb{Z}, \  \xi \in \mathcal{O}.
$$
According to the assumptions (A4) and (A5), we have  
$$ | K^x\left( t,\xi \right) |\le L_{f_2}<\infty,\quad | H^x\left( t,\xi ,z \right) |\le L_{g_2}<\infty,\quad J^x( t,\xi ) \ge 0. $$ 

Now, we introduce the following auxiliary problem 
\begin{eqnarray}\label{or418}
\begin{split}
\begin{cases}
d\vartheta \left( t \right) &=\left( A_2\left( t \right) -\alpha \right) \vartheta  \left( t \right) dt+K^x\left( t \right) \vartheta  \left( t \right) d\bar{W}^{Q_2}\left( t \right)   +\int_{\mathbb{Z}}{H^x\left( t,z \right) \vartheta  \left( t \right)}\tilde{N}_{{2}^{'}}\left( dt,dz \right)  \\
\vartheta  \left( s \right) &=y_s,
\end{cases} 
\end{split}
\end{eqnarray}
where $ y_s\in L^p\left( \varOmega ;E \right) $ is $ \mathcal{F}_s $-measurable. Denote the solution of (\ref{or418}) by  $ \vartheta \left( t;s,y_s \right) $ 
\begin{eqnarray}
\vartheta  \left( t;s,y_s \right) &=&U_{\alpha }\left( t,s \right) y_s+\psi _{\alpha }\left( \vartheta  \left( \cdot;s,y_s \right) ;s \right) \left( t \right) +\int_s^t{U_{\alpha }\left( t,r \right) K^x\left( r \right) \vartheta  \left( r;s,y_s \right)}d\bar{W}^{Q_2}\left( r \right) \cr
&&+\int_s^t{\int_{\mathbb{Z}}{U_{\alpha }\left( t,r \right) H^x\left( r,z \right) \vartheta  \left( r;s,y_s \right)}}\tilde{N}_{{2}^{'}}\left( dr,dz \right). \nonumber
\end{eqnarray}
Take expectation of the above equation and multiply both two sides  by $ e^{\delta p\left( t-s \right)} $.  Because $ \alpha $ is large enough, according to the Lemma 2.4 in \cite{cerrai2017averaging},  we have
\begin{align}
e^{\delta p\left( t-s \right)}\mathbb{E} \left\| \vartheta \left( t;s,y_s \right) \right\|^{p}
&\leq C_pe^{-\left( \alpha-\delta\right)  p\left( t-s \right)}\mathbb{E}\big\| e^{\gamma_2 \left( t,s \right)A_2}y_s \big\|^{p}+C_pe^{ \delta p\left( t-s \right)}\mathbb{E}\big\|\psi _{\alpha }\left( \vartheta  \left( \cdot;s,y_s \right) ;s \right) \left( t \right)  \big\|^{p}\cr
&\quad+C_pL_{f_2}^{p}\mathbb{E}\Big\|  \int_s^t{e^{\gamma_2 \left( t,r \right)A_2}e^{-\left( \alpha-\delta\right) \left( t-r \right)}e^{\delta\left( r-s\right) }\vartheta \left( r;s,y_s \right)}d\bar{W}^{Q_2}\left( r \right) \Big\|^{p}\cr
&\quad+C_pL_{g_2}^{p}\mathbb{E}\Big\| \int_s^t{\int_{\mathbb{Z}}{e^{\gamma_2 \left( t,r \right)A_2}e^{-\left( \alpha-\delta\right) \left( t-r \right)}e^{\delta\left( r-s\right) }\vartheta \left( r;s,y_s \right)}}\tilde{N}_{{2}^{'}}\left( dr,dz \right) \Big\|^{p}\cr
&\leq C_pL_{f_2}^{p}\underset{r\in \left[ s,t \right]}{\sup}e^{\delta p \left( r-s \right)}\mathbb{E}\left\| \vartheta \left( r;s,y_s \right) \right\|^{p}
\Big( \int_s^t{\big\| e^{\gamma_2 \left( t,r \right)A_2}e^{\left( \delta -\alpha \right) \left( t-r \right)}Q_2 \big\| _{2}^{2}}dr \Big)^{\frac{p}{2}}\cr
&\quad +C_pL_{g_2}^{p}\underset{r\in \left[ s,t \right]}{\sup}e^{\delta p \left( r-s \right)}\mathbb{E}\left\| \vartheta \left( r;s,y_s \right) \right\|^{p}
\Big( \int_s^t{\int_{\mathbb{Z}}{\big\| e^{\gamma_2 \left( t ,r\right)A_2}e^{\left( \delta -\alpha \right) \left( t-r \right)} \big\|^{2}}} v_{{2}^{'}}\left( dz \right) dr \Big)^{\frac{p}{2}}\cr
&\quad+C_pL_{g_2}^{p}\underset{r\in \left[ s,t \right]}{\sup}e^{\delta p \left( r-s \right)}\mathbb{E}\left\| \vartheta \left( r;s,y_s \right) \right\|^{p}
\int_s^t{\int_{\mathbb{Z}}{\big\| e^{\gamma_2 \left( t,r \right)A_2}e^{\left( \delta -\alpha \right) \left( t-r \right)} \big\|^{p}}} v_{{2}^{'}}\left( dz \right) dr \cr
&\quad +C_p \big\|  y_s \big\|^{p}:=\sum_{i=1}^3{\mathcal{I}_{i} \left( t\right) }+C_p \big\|  y_s \big\|^{p}.\nonumber
\end{align}
Then, the argument used in the proof of Lemma 4.2 in \cite{Xu2018Averaging} can be adapted to the present situation, and it is possible to show that if we take  $ \bar{p}>1 $ such that  $ \frac{\beta _2\left( \rho _2-2 \right)}{\rho _2}\frac{\bar{p}}{\bar{p}-2}<1, $ for any  $p\ge \bar{p}$ and $ 0<\delta<\alpha $,  we can get
$$
\sum_{i=1}^3{\mathcal{I}_{i} \left( t\right)}\leq C_{p,1}\frac{L^p}{\left( \alpha -\delta \right) ^{C_{p,2}}}\underset{r\in \left[ s,t \right]}{\sup}e^{\delta p\left( r-s \right)}\mathbb{E}\left\| \vartheta \left( r;s,y_s \right) \right\|^{p},
$$
where $ L=\max \left\{ L_{f_2},L_{g_2} \right\}. $ Hence 
$$
\underset{r\in \left[ s,t \right]}{\sup}e^{\delta p\left( r-s \right)}\mathbb{E}\left\| \vartheta \left( r;s,y_s \right) \right\| ^{p}\le C_p\left\| y_s \right\| ^{p}+C_{p,1}\frac{L^p}{\left( \alpha -\delta \right) ^{C_{p,2}}}\underset{r\in \left[ s,t \right]}{\sup}e^{\delta p\left( r-s \right)}\mathbb{E}\left\| \vartheta \left( r;s,y_s \right) \right\|^{p}.
$$
For $ \alpha >0 $ large enough, we can find  $ 0<\bar{\delta}_p<\alpha  $, such that
$$
C_{p,1}\frac{L^p}{\left( \alpha -\bar{\delta}_p \right) ^{C_{p,2}}}<1.
$$
This implies that
$$
\underset{r\in \left[ s,t \right]}{\sup}e^{p\bar{\delta}_p\left( r-s \right)}\mathbb{E}\left\| \vartheta \left( r;s,y_s \right) \right\| ^{p}\le C_p\left\| y_s \right\|^{p}.
$$
Then, let $ \delta _p=p\bar{\delta}_p $, we have
\begin{eqnarray}\label{en419}
\mathbb{E}\left\| \vartheta \left( r;s,y_s \right) \right\| ^{p}\le C_pe^{-\delta _p\left( r-s \right)}\left\| y_s \right\| ^{p}, \quad s<r.
\end{eqnarray}

Next, for any $ \mathcal{F}_s $-measurable $ y_s\in L^p\left( \varOmega ;E \right), $ we introduce the following equation 
\begin{eqnarray}\label{or417}
\begin{split}
\begin{cases}
d\varrho \left( t \right) &=\left[ \left( A_2\left( t \right) -\alpha \right) \varrho \left( t \right) -J^x\left( t \right) \varrho \left( t \right) \right] dt+K^x\left( t \right) \varrho \left( t \right) d\bar{W}^{Q_2}\left( t \right)\\ 
&\quad +\int_{\mathbb{Z}}{H^x\left( t,z \right) \varrho \left( t \right)}\tilde{N}_{{2}^{'}}\left( dt,dz \right),\\
\varrho \left( s \right) &=y_s,
\end{cases} 
\end{split}
\end{eqnarray}
denote its solution by $ \varrho \left( t;s,y_s \right). $ 
Due to the equation (\ref{or417}) is linear and $ J^x\left(t,\xi \right)\geq 0, $  using the comparison argument \cite{donati1993white} for it, we can get 
\begin{align}
y_s \ge 0,\quad  \mathbb{P}-a.s. \Rightarrow 0\leq \varrho \left( t;s,y_s \right) \leq \vartheta \left( t;s,y_s \right) ,\quad s<t,\ \mathbb{P}-a.s. 
\end{align} 
Moreove, due to the linearity of the equation (\ref{or417}), we can conclude
\begin{align}\label{en420}
&\quad Y^x\left( t;s,y \right) -Y^x\left( t;s-h,y \right) =\varrho \left( t;s,y-Y^x\left( s;s-h,y \right) \right) \cr
&=\varrho \left( t;s,y-Y^x\left( s;s-h,y \right) \land y \right) -\varrho \left( t;s,Y^x\left( s;s-h,y \right) -Y^x\left( s;s-h,y \right) \land y \right). 
\end{align}
Then, thanks to (\ref{en419})-(\ref{en420}) and (\ref{en413}),  we can get that there exists some $ \delta_p>0, $ such that
\begin{align}\label{en417}
\mathbb{E}\left\|  Y^x\left( t;s,y \right) -Y^x\left( t;s-h,y \right) \right\| ^{p}
&\le \mathbb{E}\left\|  \varrho \left( t;s,y-Y^x\left( s;s-h,y \right) \land y \right) \right\| ^{p} \cr
&\quad+\mathbb{E}\left\| \varrho \left( t;s,Y^x\left( s;s-h,y \right) -Y^x\left( s;s-h,y \right) \land y \right) \right\| ^{p} \cr
&\le \mathbb{E}\left\|  \vartheta \left( t;s,y-Y^x\left( s;s-h,y \right) \land y \right) \right\| ^{p} \cr
&\quad+\mathbb{E}\left\| \vartheta \left( t;s,Y^x\left( s;s-h,y \right) -Y^x\left( s;s-h,y \right) \land y \right) \right\| ^{p} \cr
&\le C_pe^{-\delta _p\left( t-s \right)}\mathbb{E}\lVert y-Y^x\left( s;s-h,y \right) \rVert ^p\cr
&\le C_pe^{-\delta _p\left( t-s \right)}\big( 1+\lVert x \rVert ^p+\lVert y \rVert ^p+e^{-\delta ph}\lVert y \rVert ^p \big).
\end{align}
Thanks to the completeness $ L^p\left( \varOmega ;E \right), $ it allows us to conclude that if we let $ s\rightarrow -\infty, $ there exists some $ \eta ^x\left( t \right) \in L^p\left( \varOmega ;E \right) $  such that (\ref{en414}) hold. Moreover, if we let  $ h\rightarrow \infty  $ in (\ref{en417}), we can get (\ref{en415}). 

Finally, using the same arguments as the Lemma 4.2 in our previous work \cite{Xu2018Averaging}, it is possible to prove that the limit $ \eta ^x\left( t \right) $  does not depend on the initial condition and $ \eta ^x\left( t \right) $  is a mild solution of (\ref{en42}). Moreover, by proceeding as the proof of Proposition 5.4 in \cite{cerrai2017averaging} and (\ref{en417}) in this paper, we can also study the dependence of $ \eta^x $ on the parameter $ x\in E  $ and get the conclusion (\ref{419}). We will not give specific proof here. \qed

\bibliography{references}

\begin{thebibliography}{10}
\expandafter\ifx\csname url\endcsname\relax
  \def\url#1{\texttt{#1}}\fi
\expandafter\ifx\csname urlprefix\endcsname\relax\def\urlprefix{URL }\fi
\expandafter\ifx\csname href\endcsname\relax
  \def\href#1#2{#2} \def\path#1{#1}\fi

\bibitem{Bogolyubov1961asymptotic}
N.~N. Bogolyubov, Y.~A. Mitropolskii, Asymptotic Methods in the Theory of
  Nonlinear Oscillations, Gordon and Breach Science Publishers, New York, 1961.

\bibitem{khas1968on}
R.~Khasminskii, On the averaging principle for stochastic differential it\^{o}
  equations, Kybernetika. 4 (1968) 260--279.

\bibitem{givon2007strong}
D.~Givon, Strong convergence rate for two-time-scale jump-diffusion stochastic
  differential systems, Multiscale. Model. Sim. 6 (2007) 577--594.

\bibitem{freidlin2012random}
M.~Freidlin, A.~Wentzell, Random Perturbations of Dynamical Systems, Springer
  Science and Business Media, Berlin Heidelberg, 2012.

\bibitem{duan2014effective}
J.~Q. Duan, W.~Wang, Effective Dynamics of Stochastic Partial Differential
  Equations, Elsevier, 2014.

\bibitem{xu2011averaging}
Y.~Xu, J.~Q. Duan, W.~Xu, An averaging principle for stochastic dynamical
  systems with {L\'{e}vy} noise, Physica D. 240 (2011) 1395--1401.

\bibitem{xu2015approximation}
Y.~Xu, B.~Pei, Y.~G. Li, Approximation properties for solutions to
  {non-Lipschitz} stochastic differential equations with {L\'{e}vy} noise,
  Math. Method Appl. Sci. 38 (2015) 2120--2131.

\bibitem{xu2017stochastic}
Y.~Xu, B.~Pei, J.~L. Wu, Stochastic averaging principle for differential
  equations with {non-Lipschitz} coefficients driven by fractional {Brownian}
  motion, Stoch. Dynam. 17 (2017) 1750013.

\bibitem{Cerrai2009Averaging}
S.~Cerrai, M.~Freidlin, Averaging principle for a class of stochastic
  reaction-diffusion equations, Probab. Theory Relat. Fields. 144 (2009)
  137--177.

\bibitem{cerrai2009khasminskii}
S.~Cerrai, A khasminskii type averaging principle for stochastic
  reaction-diffusion equations, Ann. Appl. Probab. 19 (2009) 899--948.

\bibitem{cerrai2011averaging}
S.~Cerrai, Averaging principle for systems of reaction-diffusion equations with
  polynomial nonlinearities perturbed by multiplicative noise, SIAM. J. Math.
  Anal. 43 (2011) 2482--2518.

\bibitem{wang2012average}
W.~Wang, A.~Roberts, Average and deviation for slow-fast stochastic partial
  differential equations, J. Differ. Equations 253 (2012) 1265--1286.

\bibitem{Fu2011AN}
H.~FU, J.~Duan, An averaging principle for two-scale stochastic partial
  differential equations, Stoch. Dyn. 11~(2-3) (2011) 353--367.

\bibitem{Bao2017Two}
B.~J., Y.~G., C.~Yuan, Two-time-scale stochastic partial differential equations
  driven by $ \alpha $-stable noises: averaging principles, Bernoulli 23 (2017)
  645–669.

\bibitem{dong2018averaging}
Z.~Dong, X.~Sun, H.~Xiao, J.~Zhai, Averaging principle for one dimensional
  stochastic {Burgers} equation, J. Differ. Equations 265~(10) (2018)
  4749--4797.

\bibitem{Fu2015Strong}
H.~Fu, L.~Wan, J.~Liu, Strong convergence in averaging principle for stochastic
  hyperbolic–parabolic equations with two time-scales, Stoch. Process. Appl.
  125~(8) (2015) 3255--3279.

\bibitem{pei2017two}
B.~Pei, Y.~Xu, J.~L. Wu, Two-time-scales hyperbolic-parabolic equations driven
  by {Poisson} random measures: Existence, uniqueness and averaging principles,
  J. Math. Anal. Appl. 447 (2017) 243--268.

\bibitem{gao2019averaging}
P.~Gao, Averaging principle for multiscale stochastic {Klein}-{Gordon}-{Heat}
  system, J. Nonlinear Sci. 29~(4) (2019) 1701--1759.

\bibitem{gao2018averaging}
P.~Gao, Averaging principle for the higher order nonlinear {Schr\"{o}dinger}
  equation with a random fast oscillation, J. Stat. Phys. 171~(5) (2018)
  897--926.

\bibitem{cerrai2017averaging}
S.~Cerrai, A.~Lunardi, Averaging principle for non-autonomous slow-fast systems
  of stochastic reaction-diffusion equations: The almost periodic case, SIAM.
  J. Math. Anal. 49 (2017) 2843--2884.

\bibitem{liu2020averaging}
W.~Liu, M.~R{\"o}ckner, X.~Sun, Y.~Xie, Averaging principle for slow-fast
  stochastic differential equations with time dependent locally lipschitz
  coefficients, J. Differ. Equations 268~(6) (2020) 2910--2948.

\bibitem{Xu2018Averaging}
Y.~Xu, R.~F. Wang, Averaging principles for non-autonomous two-time-scale
  stochastic reaction-diffusion equations with jump, Complexity (2020) In
  Press.

\bibitem{prato2014stochastic}
G.~D. Prato, J.~Zabczyk, Stochastic Equations in Infinite Dimensions, Cambridge
  University Press, Cambridge, 2014.

\bibitem{cerrai2003stochastic}
S.~Cerrai, Stochastic reaction-diffusion systems with multiplicative noise and
  non-{Lipschitz} reaction term, Probab. Theory Relat. Fields. 125 (2003)
  271--304.

\bibitem{peszat2007stochastic}
S.~Peszat, J.~Zabczyk, Stochastic Partial Differential Equations with L\'{e}vy
  noise: An Evolution Equation Approach, Cambridge University Press, Cambridge,
  2007.

\bibitem{Applebaum2009Processes}
D.~Applebaum, L\'{e}vy Processes and Stochastic Calculus, Cambridge University
  Press, Cambridge, 2009.

\bibitem{Pazy2012Semigroups}
A.~Pazy, Semigroups of linear operators and applications to partial
  differential equations, Springer Science \& Business Media, 2012.

\bibitem{Fu2011Strong}
H.~Fu, J.~Liu, Strong convergence in stochastic averaging principle for two
  time-scales stochastic partial differential equations, J. Math. Anal. Appl.
  384~(1) (2011) 70--86.

\bibitem{donati1993white}
C.~Donati-Martin, E.~Pardoux, White noise driven spdes with reflection, Probab.
  Theory Rel. 95 (1993) 1--24.

\end{thebibliography}
\end{document}